\scrollmode
\documentclass[12pt]{amsart}
\usepackage{verbatim}
\usepackage{eucal}
\usepackage{amssymb}
\usepackage{graphicx}
\usepackage{psfrag}
\usepackage{mathrsfs}
\usepackage[all]{xy}
\newif \ifwide
\newif \ifavnermargin
\def \makemargins{
\ifwide
	\oddsidemargin .25in
	\evensidemargin .25in
	\textwidth 6.00in
\else
\fi
\ifavnermargin
	\headheight=7pt
	\textheight=574pt
	\textwidth=432pt
	\topmargin=14pt
	\oddsidemargin=18pt
	\evensidemargin=18pt
\else	
\fi
}

\avnermargintrue
\makemargins
\theoremstyle{plain}
\newtheorem{proposition}[subsection]{Proposition}
\newtheorem{lemma}[subsection]{Lemma}

\theoremstyle{definition}
\newtheorem{definition}[subsection]{Definition}

\theoremstyle{remark}
\newtheorem{remark}[subsection]{Remark}

\newcount\timehh\newcount\timemm
\timehh=\time 
\divide\timehh by 60 \timemm=\time
\newcommand{\draftauthor}[1]{\author{#1
    {
      --- \protect \protect\sc\today\ ---
      \ifnum\timehh<10 0\fi\number\timehh\,:\,\ifnum\timemm<10 0\fi\number\timemm
      \protect \, \, \protect \bf DRAFT
    }
  }
}
\newcommand{\R}{{\mathbb R}}
\newcommand{\F}{{\mathbb F}}
\newcommand{\C}{{\mathbb C}}
\newcommand{\Z}{{\mathbb Z}}
\newcommand{\z}{{\mathbb Z}}
\newcommand{\Q}{{\mathbb Q}}

\newcommand{\V}{{\mathscr V}}
\newcommand{\Proj}{{\mathbb P}}
\renewcommand{\P}{{\mathbb P}}
\newcommand{\PPP}{\mathscr P}

\newcommand{\aaa}{\mathbf{a}}
\renewcommand{\H}{\mathscr H}

\def\bs{\backslash}
\def\op{\operatorname}
\DeclareMathOperator{\SL}{SL}
\DeclareMathOperator{\GL}{GL}
\DeclareMathOperator{\SO}{SO}

\def\vor{Vorono\v{\i}}

\def\scrV{\V}
\def\scrP{\PPP}

\newcommand{\Vor}{Vorono\v{\i}}

\newcommand{\vv}{{\mathbf {v}}}
\renewcommand{\aa}{{\mathbf {a}}}

\newcommand{\uu}{{\mathbf {u}}}
\newcommand{\yy}{{\mathbf {y}}}
\newcommand{\qq}{{\mathbf {q}}}

\newcommand{\coinv}[1]{(S_{#1})_{\Gamma }}

\DeclareMathOperator*{\Max}{Max}

\DeclareMathOperator{\support}{supp}

\DeclareMathOperator{\Frob}{Frob}

\DeclareMathOperator{\sign}{sgn}

\DeclareMathOperator{\Gal}{Gal}
\CompileMatrices

\begin{document}

\title{Cohomology of congruence subgroups of $\SL_{4}(\Z)$}

\newif \ifdraft
\def \makeauthor{
\ifdraft
	\draftauthor{Avner Ash, Paul Gunnells, Mark McConnell}
\else

\author{Avner Ash}
\address{Department of Mathematics\\
The Ohio State University\\
231 W. 18th Ave\\
Columbus, OH 43210}
\email{ash@math.ohio-state.edu}

\author{Paul E. Gunnells}
\address{Department of Mathematics\\
Columbia University\\
New York, NY 10027}
\email{gunnells@math.columbia.edu}

\author{Mark McConnell}
\address{WANDL, Inc.\\
2121 Route 22 West\\
Bound Brook, NJ  08805}
\email{mmcconnell@wandl.com}

\fi
}

\thanks{The research of the first author was partially supported by
NSF grant DMS--9531675.  The second author was partially supported by
a Columbia University Faculty Research grant and NSF grant
DMS--9627870.  The third author was partially supported by NSF grant
DMS--9704535.}
\draftfalse
\makeauthor

\ifdraft
	\date{\today}
\else
	\date{March 30, 2000}
\fi

%
%
\begin{abstract}
Let $N>1$ be an integer, and let $\Gamma = \Gamma _{0} (N) \subset
\SL_{4} (\Z)$ be the subgroup of matrices with bottom row congruent to
$(0,0,0,*)\mod N$.  We compute $H^{5} (\Gamma; \C) $ for a range of
$N$, and compute the action of some Hecke operators on many of these
groups.  We relate the classes we find to classes coming from the
boundary of the Borel-Serre compactification, to Eisenstein series,
and to classical holomorphic modular forms of weights 2 and 4.

\end{abstract}
\maketitle

%
%
\section{Introduction}\label{introduction}
\subsection{}\label{diag.sect}
Let $n\ge 1$, and let $\Gamma$ be a congruence subgroup of $\SL_n(\z)$
of level $N$.  Let $S_N$ be the subsemigroup of the integral matrices
in $\GL_{n}(\Q)$ such that $(\Gamma, S_N)$ is a Hecke pair.

We denote by $\H(N)$ the $\C$-algebra of double cosets $\Gamma S_N
\Gamma$.  This algebra acts on the cohomology and homology of $\Gamma$
with any coefficient $\C S_N$-module.  When a double coset is acting,
we call the map defined by its action a Hecke operator.  In this paper
we will work only with the trivial coefficient module $\C$.

Let $l$ be a prime not dividing $N$, and let $D (l,k)$ be the diagonal
matrix with diagonal $(1,\dots ,1,l,\dots ,l)$, where the number of
$l$'s is $k$.  Then $\H(N)$ contains all double cosets of the form
$\Gamma D(l,k) \Gamma$, and we denote the corresponding Hecke operator by
$T (l,k)$.  Fix a prime $p$ not dividing $N$, and an embedding of $\Q
_{p}$ into $\C$.  Let $G_{\Q } = \Gal (\bar \Q /\Q )$.

\begin{definition}\label{attached.def}
Let $\mathcal V$ be an $\H(N)$-module and suppose $\beta \in \mathcal
V$ is an eigenclass for the action of $\H(N)$.  For $l$ prime to
$N$, write $T(l,k)(\beta) =
a(l,k)\beta $ , where the $a (l,k)\in \C,
k=0,\dots ,n$ are algebraic integers.  Let $\rho$ be a continuous
semisimple representation $\rho:G_{\Q}\rightarrow \GL_{n}(\Q_p)$,
unramified outside $pN$, such that
\begin{equation}\label{star}
\sum_k (-1)^k l^{k(k-1)/2} a(l, k) X^k=\det (I-\rho(\Frob_l) X)
\end{equation}
for all $l$ not dividing $pN$.  Then we shall say that $\rho$ is
\emph{attached} to $\beta $.
\end{definition}

For example, theorems of Eichler, Shimura and Deligne imply that if $n
= 2$ and $\mathcal V$ is the Hecke-module of classical holomorphic
modular cuspforms for $\SL_{2}(\z)$, then there always exists a $\rho$
attached to any Hecke eigenform $\beta \in \mathcal V$.

Standard conjectures (for example in \cite{clozel}) state that if
$\mathcal V$ is the cuspidal cohomology of $\Gamma$ with trivial
coefficients, then any Hecke eigenclass should have an attached
$p$-adic Galois representation.  In fact, one has the stronger
conjecture that there should be an attached motive
(cf. \cite{clozel}).  We can also extend the conjecture of
\cite{clozel} to include all the cohomology---in principle, the theory
of Eisenstein series should allow one to reduce the extended
conjecture to the one for cuspidal cohomology.

In a series of papers \cite{exp.ind, apt, aac, ash.tiep}, the first
author with a number of coworkers has tested this conjecture and a mod
$p$ variant of it when $n=3$.  Other tests for $n=3$ can be found in
the work of van Geemen and Top with van der Kallen and Verberkmoes
\cite{vgt1, vgt2, vgt3}.  The purpose of this paper is to make the
first computational tests of this conjecture for $n=4$.

\subsection{}
We work here with $\Gamma = \Gamma_0(N)$, defined as the subgroup of
$\SL_{n}(\z)$ consisting of the matrices with last row congruent
to $(0,\dots, 0, *)$ modulo $N$.  

Conceptually, there is no big difference in computing the cohomology of
$\Gamma$ between the cases of $n = 3$---or even $n=2$---and the case of
general $n$.  One needs to write down a simplicial complex $C$
homotopic to the chains of a $B\Gamma(1)$-space, and compute the
cohomology.

When $n=2$ or $3$, the most interesting part of cohomology, namely the
cuspidal part, occurs in the top dimension of $C$.  Therefore to
compute it all we
have to compute is the cokernel of a coboundary map.  However, when $n=4$
the cuspidal cohomology occurs in dimensions 4 and 5, whereas the
virtual cohomological dimension of $\Gamma$, and hence smallest
possible dimension of $C$, is 6.  So now the cohomology is a
subquotient, which adds considerably to the complexity of the computer
programs.  The cuspidal cohomology in degree 4 is dual to that in
degree 5, so we concentrate on computing the latter in
this paper.  

But the big difference between $n=2,3$ and $n=4$ occurs when we try to
compute the Hecke action on the cohomology.  In the top dimension we
can use the Ash-Rudolph algorithm or its variants \cite{ash.rudolph},
as was done for $n=3$ in the works cited above.  However, for $n=4$,
where we look just below the top dimension, a brand new idea was
necessary.  This is due to the second author, and is the subject of
\cite{gunn}.  Thus this paper is also a test of the algorithms
proposed in \cite{gunn}, and they pass with flying colors.  It is an
open question whether the algorithms of \cite{gunn} terminate in a
finite number of steps; in practice, though, they always terminate
quickly, and we used them here without problems.

Our computations of $H^5(\Gamma_0(N);\C)$, detailed in
\S\ref{num.results.sect}, were made for $N \le 53$.\footnote{Actually,
to avoid numerical instability in floating-point computations, we
replace $\C$ with the finite field $\F_{31991}$
(cf. \S\ref{lancz.sect}).}  No cusp forms were discovered, but some
interesting phenomenology of the boundary cohomology was observed, as
discussed in \S\ref{interpretation}.  This leads to some open
questions about the cohomology of the boundary, which are discussed
there.  Since we have not completed the Hecke computations for some
high levels near $53$, it is possible that we do have a cusp form that
we haven't yet identified as such.  There is a bound due to Fermigier
\cite{ferm} that states that if $N < 31$, there cannot be any cuspidal
cohomology.

For high levels our computations produced large sparse matrices,
as large as $110464\times 30836$ for level $48$.  To perform linear algebra
with these matrices, we used a version of the Lanczos algorithm mod
$p$, in the spirit of LaMacchia-Odlyzko \cite{ody} (see also
\cite{teit}).

\subsection{}
We thank Eric Conrad for some assistance with programming, 
and Peter Woit for excellent computing support.  We thank David
Ginzburg for conversations at the beginning of this project. 

\section{Background}\label{bkground}

\subsection{}\label{vcd.section}
Let $V$ be the $\R $-vector space of all symmetric $n\times n$
matrices, and let $C\subset V$ be the cone of positive-definite
matrices.  Then the group $G=\SL_{n} (\R )$ acts on $C$ on the left by
$(g,c)\mapsto g\cdot c\cdot g^{t}$, and the stabilizer of any given
point is isomorphic to $\SO_{n}$.

Let $X$ be $C$ mod homotheties.  The $G$-action on $C$ commutes with
the homotheties and induces a transitive $G$-action on $X$.  The
stabilizer of any given point of $X$ is again isomorphic to $\SO_{n}$.
After choosing a basepoint, we may identify $X$ with the global
riemannian symmetric space $\SL_{n} (\R )/\SO_{n}$, a contractible,
noncompact, smooth manifold of real dimension $d = n (n+1)/2 - 1$.

The group $\SL_{n} (\Z )$ acts on $X$ via the $G$-action, and does so
properly discontinuously.  Hence if $\Gamma \subset \SL_{n} (\Z )$ is
any finite-index subgroup, the quotient $\Gamma \backslash X$ is a
real noncompact manifold except for at most finitely many quotient
singularities.  We may then identify the complex group cohomology
$H^{*} (\Gamma; \C )$ with $H^{*} (\Gamma \backslash X;\C )$.
Although the dimension of $\Gamma \backslash X$ is $d$, it can be
shown that $H^{i} (\Gamma \backslash X;\C )=0$ if $i> d-n+1$
\cite[Theorem 11.4.4]{borel.serre}.  The number $\nu = d-n+1$ is
called the \emph{virtual cohomological dimension} of $\Gamma $.

In this paper we will always take $\Gamma $ to be the congruence
subgroup $\Gamma _{0} (N)$ of matrices whose last row is congruent to
$(0,\dots ,0,*) \mod N$.

\subsection{}\label{vor.poly.section}
Recall that a point in $\Z ^{n}$ is said to be \emph{primitive} if the
greatest common divisor of its coordinates is $1$.  In particular, a
primitive point is nonzero.  Let $\PPP\subset \Z ^{n}$ be the set of
primitive points.  Any $v\in \PPP$, written as a column vector,
determines a rank-one symmetric matrix $q (v)\in \bar C$ by $q ( v) =
v\cdot v^{t}$.  The \emph{\Vor \ polyhedron} $\Pi $ is the closed
convex hull of the points $q (v)$, as $v$ ranges over $\PPP$.

Note that, by construction, $\SL_{n} (\Z )$ acts on $\Pi $.  The cones
over the faces of $\Pi $ form a fan $\V$ that induces a $\Gamma
$-admissible decomposition of $C$~\cite[p. 117]{ash}.  Essentially,
this means that $\Gamma $ acts on $\V$; that each cone is spanned by a
\emph{finite} collection of points $q (v)$ where $v\in \PPP$; and that
there are only finitely many $\Gamma $-orbits in $\V$.  The fan
$\V$ provides a reduction theory for $C$ in the following sense: any
point $x\in C$ is contained in a unique $\sigma (x) \in \V$, and the
set $\{\gamma \in \SL _{n} (\Z ) \mid \gamma \cdot \sigma (x) = \sigma
(x) \}$ is finite.

\subsection{}\label{wrr.section}
We summarize facts about the well-rounded retract of \cite{Ash80,
Ash84}.  There is a deformation retraction $C\to C$ that is
equivariant under the actions of both $\bar\Gamma = \SL_n(\Z)$ and the
homotheties.  Its image modulo homotheties is the \emph{well-rounded
retract} $W$ in~$X$.  The well-rounded retract is contractible, since
it is a deformation retract of the contractible space~$X$.  Hence the
cohomology of $\Gamma\bs X$ with coefficients in $\C$ is canonically
isomorphic to the equivariant cohomology $H^i_\Gamma(W; \C)$ where
$\Gamma $ acts trivially on the coefficient module $\C$.  This is in turn
canonically isomorphic to the complex cohomology $H^{i} (\Gamma \bs
W; \C)$, since $\C$ has characteristic zero, which moreover is
isomorphic to $H^{i} (\Gamma ;\C)$.  We will focus on computation of
the equivariant cohomology.  The dimension of~$W$ equals the virtual
cohomological dimension $\nu$, and the quotient $\Gamma\bs W$ is
compact.

The well-rounded retract~$W$ is naturally a locally finite cell
complex, the cells being convex polytopes in~$V$.  The group $\SL_{n}
(\Z )$ preserves the cell structure, and the stabilizer of each cell
in $\SL _{n} (\Z )$ is finite.  The theory of cores and co-cores in
\cite[Chapter 2]{A-M-R-T} shows that the cells in $W$ are in a
one-to-one, inclusion-reversing correspondence with the cones in the
\vor{} fan $\scrV$.  By abuse of notation, the cell in $W$
corresponding to $\sigma$ will still be denoted $\sigma$.

\subsection{}
One can give a more precise description of the combinatorics of the
cones in $\scrV$ and the cells in $W$.  To each $\sigma\in \scrV$, we
define the set 
\[
M(\sigma) = \{v\in\scrP \mid \hbox{$q(v)$~is a vertex of
the face of~$\Pi$ generating $\sigma $.}\}
\]
We associate the same set $M(\sigma)$ to the corresponding cell
in~$W$, and call $M(\sigma)$ the set of \emph{minimal vectors
of\/}~$\sigma$ (because of how $W$ is constructed in
\cite{Ash80,Ash84}).  Since $\Pi$ is the convex hull of the $q(v)$'s,
it is clear that $\mu \colon \sigma \rightarrow M(\sigma)$ is an
inclusion-preserving (respectively, inclusion-reversing) bijection
between the cones in $\scrV$ (resp., cells in $W$) and a collection of
finite subsets of $\scrP$.  In principle, this reduces the study of
the combinatorics of $\scrV$ and $W$ to the study of the image of~$\mu
$.  For instance, face relations $\tau\subseteq\sigma$ in~$\scrV $ are
read off from subset relations $M(\tau)\subseteq M(\sigma)$.  In
practice, determining the image of~$\mu $ calls for explicit
computations with real quadratic forms, computations whose difficulty
grows exponentially as a function of~$n$.  The computations have been
carried out completely for $n \le 5$ by various authors.

For the rest of this subsection, we set $n=4$ and give more details.
We state results for the well-rounded retract; these imply their
analogues for $\scrV$.  The image of~$\mu $ was computed independently by
\cite{L-S} and (in essence) \cite{vSto}.  The cells of~$W$ fall into
eighteen equivalence classes modulo $\SL_{n} (\Z )$.  Let $T$ (for
``type'') be a variable running through these eighteen classes.  This
partitions the set of cells of~$W$ into eighteen pieces called the
$W_T$.  Any $\sigma\in W_T$ is said to be \emph{of type~$T$.}  In each
$W_T$, we fix one representative cell $\sigma_T$, the \emph{standard
cell of type $T$}.  The $M(\sigma_T)$'s are written down explicitly in
\cite{M91}.\footnote{Since $v\in M(\sigma) \Leftrightarrow
-v \in M(\sigma)$, it is customary to write down only one member
of the pair $\pm v$.} This determines the image of~$\mu $, since
the image is the union of the $\SL_{n} (\Z )$-translates of the eighteen
$M(\sigma_T)$'s.

\subsection{}
To compute the action of the Hecke operators on cohomology, the
well-rounded retract is insufficient, since the operators do not act
cellularly.  To ameliorate this, we use the sharbly complex.
The material in this subsection closely follows \cite{ash.sharb}.

\begin{definition}\label{sharbly.complex}
\cite{ash.sharb} The \emph{sharbly complex} is the chain complex
$\left\{S_{*},\partial \right\}$ given by the following data:
\begin{enumerate}
\item For $k\geq 0$, $S_{k}$ is the module of formal $\Z$-linear
combinations of basis elements $\uu  = [v_{1},\ldots,v_{n+k}]$, where each
$v_{i}\in \PPP$, mod the relations:
\begin{enumerate}
\item If $\tau $ is a permutation on $(n+k)$ letters, then
\[
[v_{1},\ldots,v_{n+k}] = \sign (\tau ) [\tau (v_{1}),\ldots,\tau
(v_{n+k})], 
\]
where $\sign (\tau )$ is the sign of $\tau $.
\item  If $q = \pm 1$, then 
\[
[q v_{1},v_{2}\ldots,v_{n+k}] = [v_{1},\ldots,v_{n+k}].
\]
\item If the rank of the matrix
$(v_{1},\ldots,v_{n+k})$ is less than $n$, then $\uu = 0$.
\end{enumerate}
\item The boundary map $\partial \colon S_{k}\rightarrow S_{k-1}$ is 
$$[v_{1},\ldots,v_{n+k}] \longmapsto \sum _{i=1}^{n+k}  (-1)^{i}
[v_{1},\ldots,\hat{v_{i}},\ldots,v_{n+k}].  $$
\end{enumerate}
\end{definition}

The basis elements $\uu = [v_{1},\dots ,v_{n+k}]$ are called
\emph{$k$-sharblies}.  By abuse of notation, we will often use the
same symbol $\uu $ to denote a $k$-sharbly and the $k$-sharbly chain
$1\cdot \uu $.  The obvious left action of $\Gamma $ on $S_{* }$
commutes with $\partial $.

For any $k\geq 0$, let $\coinv{k}$ be the module of $\Gamma
$-coinvariants.  This is the quotient of $S_{k}$ by the relations of
the form $\gamma \cdot \uu - \uu $, where $\gamma \in \Gamma $, $\uu
\in S_{k}$.  This is also a complex with the induced boundary, which
we denote by $\partial_{\Gamma } $.  It is known (cf. \cite{gunn})
that $H^{\nu - k}(\Gamma ; \C)$ is naturally isomorphic to $H_{k}
(\coinv{*}\otimes \C )$.

Let $\uu = [v_{1},\dots ,v_{n+k}]$ be a $k$-sharbly.  Let $\|\uu
\|$ be 
\[
\Max |\det (v_{i_{1}}, \dots , v_{i_{n}})|,
\]
where the maximum is taken over all $n$-fold subsets $\{i_{1},\dots
,i_{n} \}\subset \{1,\dots ,n+k \}$.  Note that this quantity is
well-defined mod the relations in Definition \ref{sharbly.complex}.
We extend this notion to sharbly chains $\xi = \sum n (\uu)\uu $ by
setting $\|\xi \|$ to be the maximum of $\|\uu \|$, as $\uu $ ranges
over all sharblies in the support of $\xi $.  We say that $\xi $ is
\emph{reduced} if $\|\xi \| = 1$.  It is known (cf. \cite{M91}) that
for $\Gamma \subset \SL _{4} (\Z )$, the group $H^{5} (\Gamma ; \C)$
is spanned by reduced $1$-sharbly cycles.

\subsection{}\label{that.section}
Since the generators of the sharbly complex are indexed by sets of
primitive vectors, it is clear that there is a close relationship
between $S_{*}$ and the chain complex associated to $W$, although of
course $S_{*}$ is much bigger.  Both complexes compute $H^{*} (\Gamma
; \C)$.  We refer to \cite{gunn} for a discussion of this, phrased in
terms of the fan $\V $.  The main advantage of $\coinv{*}$ is that it
admits a Hecke action.  Specifically, let $\xi = \sum n (\uu )\uu $ be
a sharbly cycle mod $\Gamma $, and consider the Hecke operator $T
(l,k)$ associated to the double coset $\Gamma D (l,k)\Gamma $
(cf. \S\ref{diag.sect}).  Write
\[
\Gamma D (l,k)\Gamma = \coprod _{g\in \Omega } \Gamma g, 
\] 
a finite (disjoint) union.  Then 
\begin{equation}\label{how.hecke.acts}
T (l,k) (\xi ) = \sum _{g\in \Omega , \uu } n (\uu )g\cdot \uu . 
\end{equation}
Since $\Omega \not \subset \SL _{n} (\Z )$ in general, the Hecke-image of a
reduced sharbly isn't usually reduced.

\section{Implementation details}\label{imp.details.sect}
\subsection{}\label{subsectone} 
We state our results for general~$n$
as much as possible, though our main case of interest is $n=4$.  We
have working programs for $n\le4$.  Though we focus on $\SL_n(\Z)$,
analogous results hold for $\GL_n(\Z)$, and we have working programs
for both $\SL$ and $\GL$.

Section \ref{imp.details.sect} is very technical.  The reader may wish
to skip to \S\ref{heckeop.sect} or \S\ref{num.results.sect}.

\subsection{}\label{subsecttwo}
Let $\bar\Gamma = \SL _{n} (\Z )$.  Recall that $W_T$ is the
$\bar\Gamma$-orbit of cells of type $T$ in~$W$.  Let $\sigma_T$ be
a fixed representative cell in~$W_T$.  The stabilizer in $\bar\Gamma$
of $\sigma_T$ is denoted $\bar\Gamma_{\sigma_T}$, or $\bar\Gamma_T$
for short.  This is a finite group that is straightforward to compute,
since the minimal vectors $M(\sigma_T)$ are known.  Our program
maintains a database of the $\bar\Gamma_T$.

Standard facts about stabilizers give the following:

\begin{proposition}\label{propone}
There is a one-to-one correspondence between cells $\sigma \in W_T$
and cosets $\bar\Gamma / \bar\Gamma_T$, given by $\gamma\sigma_T
\leftrightarrow \gamma \bar\Gamma_T$ for any $\gamma\in \bar\Gamma$
such that $\sigma = \gamma\sigma_T$.
\end{proposition}

Under the smaller group $\Gamma$, the $\bar\Gamma$-orbit $W_T$ breaks
up into suborbits.  If $\Gamma$ were a torsion-free group, $\Gamma\bs
W$ would be a finite cell complex, its cells would be given exactly by
the $\Gamma$-suborbits, and we could compute $H^i(\Gamma\bs W)$ by the
standard methods for cell complexes.  In our case, $\Gamma$ is not
torsion-free, but $\Gamma\bs W$ can be thought of as a finite
``orbifold cell'' complex, whose elements are the quotients of cells
by finite groups; the $\Gamma$-suborbits are in one-to-one
correspondence with the orbifold cells.

The goal of this subsection is to understand the $\Gamma$-suborbits in
terms of the actions of the $\bar\Gamma_T$ on finite projective
spaces.  By $\P^{n-1} = \P^{n-1}(\Z/N\Z)$, we mean the set of vectors
$(x_1, \dots, x_n)\in (\Z/N\Z)^n$ that are primitive in the sense
that the ideal $(x_1,\dots,x_n)$ in $\Z/N\Z$ is~$(1)$, modulo the
equivalence relation given by scalar multiplication by the units
$(\Z/N\Z)^\times$ of $\Z/N\Z$.  When $N$ is a prime, this is the usual
projective space over the field of $N$ elements.  As usual, the
equivalence class of the vector $(x_1,\dots,x_n)$ is denoted $\aaa =
[x_1:\cdots:x_n]$.  We view these $n$-tuples as rows rather than
columns; $\bar\Gamma$ acts on the right on $\P^{n-1}(\Z/N\Z)$ in the
obvious way.

We define the {\it bottom row map\/} $\mathfrak{b}\colon \bar\Gamma \to
\P^{n-1}(\Z/N\Z)$ as follows.  For a matrix $\gamma\in \bar\Gamma$,
the bottom row of~$\gamma$ is a primitive vector in $\Z^n$.  Let
$\mathfrak{b}(\gamma)$ be the equivalence class of this image in
$\P^{n-1}(\Z/N\Z)$.

\begin{lemma}\label{lemone}
The bottom row map $\mathfrak{b}:\bar\Gamma \to
\P^{n-1}(\Z/N\Z)$ induces a bijection between $\Gamma\bs \bar\Gamma 
$ and $\P^{n-1}(\Z/N\Z)$, given by
$$
\Gamma\gamma \mapsto \mathfrak{b}(\gamma).
$$
The map is equivariant for the right action of $\bar\Gamma$.
\end{lemma}

\begin{proof}
It is a standard fact that a vector in $\Z^n$ is primitive if and
only if it is the bottom row of some element of $\bar\Gamma$.  This
implies $\mathfrak{b}$ is surjective and that
$\mathfrak{b}^{-1}([0:\cdots:0:0:1]) = \Gamma$.  The rest is clear.
\end{proof}

We can now describe the $\Gamma$-orbits of cells in each $W_T$.

\begin{proposition}\label{proptwo}
The $\Gamma$-orbits of cells in~$W_T$ are
in one-to-one correspondence with the orbits~$O$ of the right
$\bar\Gamma_T$-action on $\P^{n-1}$. 
\end{proposition}

\begin{proof}
 We have $\Gamma \bs W_T = \Gamma \bs \bar\Gamma / \bar\Gamma_T $ by
Proposition \ref{propone}, and this equals $\P^{n-1}/\bar\Gamma_T$ by
Lemma \ref{lemone}.
\end{proof}

The first step of our computer program is to determine, for each
type~$T$, the decomposition of $\P^{n-1}$ into right
$\bar\Gamma_T$-orbits.  Since we are primarily studying $H^i_\Gamma(W;
\C)$ for $i=5,6$, we only need to work with the~$T$ representing cells
of dimensions~4, 5 and~6.

We note the following fact, whose proof is immediate.

\begin{lemma}\label{lemtwo}
Let $\gamma_0\in\bar\Gamma$, and let $\aaa =
\mathfrak{b}(\gamma_0)$ in $\P^{n-1}$.  Then the stabilizer of $\aaa$
under the right action of $\bar\Gamma$ is $\gamma_0^{-1} \Gamma
\gamma_0$.  
\end{lemma}

Let $\sigma\in W$ be a cell of type~$T$, with $\sigma = \gamma_0
\sigma_T$.  Its stabilizer in $\Gamma$, denoted $\Gamma_\sigma$, is
clearly
\begin{equation}
\Gamma_\sigma = (\gamma_0 \bar\Gamma_T \gamma_0^{-1}) \cap \Gamma.
\end{equation}
Lemma~\ref{lemtwo} implies

\begin{lemma}\label{lemthree}
The group $\gamma_0^{-1} \Gamma_\sigma \gamma_0$ is the subgroup of
$\bar\Gamma_T$ that preserves~$\aaa=\mathfrak{b} (\gamma _{0})$ under the right
action on $\P^{n-1}$.
\end{lemma}

\subsection{}\label{subsectthree}
 In this subsection, we fix orientations on the cells $\sigma\in W$.
It's necessary to be extremely careful---mistakes in orientation are
easy to make and will ruin the computations.  The price to pay is to
sort through the details of the action of $\bar\Gamma$.

Recall that $\bar\Gamma_T$ is the stabilizer of $\sigma_T$ in
$\bar\Gamma$.  For each~$T$, there is an orientation character
$\bar\Gamma_T \to \{\pm1\}$ indicating whether or not $\gamma\in
\bar\Gamma_T$ preserves the orientation on $\sigma_T$.  Our program
stores the values of these characters along with $\bar\Gamma_T$.  Let
$\bar\Gamma_T^+$ be the subgroup of $\bar\Gamma_T$ where the
orientation is~$+1$.

\begin{remark}\label{remone}
To compute the value of the character at~$\gamma$, we determine
(i)~how $\gamma$ acts on the orientation of the cone~$C$, and divide
by (ii)~how $\gamma$ acts on the orientation of the \vor{} cone dual
to $\sigma_T$.  Dividing works because $C$ is locally the direct
product of the cell and its dual \vor{} cone.  As for~(i), every
element of $\bar\Gamma$ acts by $+1$ on the orientation of~$C$, since
$\bar\Gamma$ is a subgroup of the connected group $\SL_n(\R)$.  In
this paper, where $n=4$ and $\dim\sigma > 0$, it turns out that all
the dual \vor{} cones are simplicial; the sign in~(ii) is the sign of
the permutation that~$\gamma$ effects on the bounding rays $q(v)$ of
the cone, which is easily computed.
\end{remark}

Let $O$ be a right $\bar\Gamma_T$-orbit in $\P^{n-1}$.  We call~$O$
{\it non-orientable\/} if for some (which implies every) $\aaa\in
O$, there exists $\gamma\in \bar\Gamma_T \smallsetminus \bar\Gamma_T^+$ with $\aaa\gamma = \aaa$.  Otherwise, we call~$O$ {\it orientable}.  These
notions depend on~$T$, though we usually leave~$T$ out of the
notation.

If~$O$ is orientable, fix some $\aaa_0\in O$.  Define the {\it
orientation number\/} of $\aaa\in O$ to be $+1$ (resp., $-1$)
according as $\aaa = \aaa_0 \gamma$ for some $\gamma\in
\bar\Gamma_T^+$ (resp., $\gamma\in\bar\Gamma_T
\smallsetminus\bar\Gamma_T^+$).  The orientation number is
well-defined precisely because~$O$ is orientable.  Again, the notions
depend on the choice of $\aaa_0$, though we leave $\aaa_0$ out of the
notation.

Let $\gamma\in \bar\Gamma_T$, and let $\rho$ be any cell of~$W$ with
any given orientation.  Since
$\gamma$ acts by diffeomorphisms on $C$, it carries the orientation
on $\rho$ to some orientation on the cell $\gamma\rho$.  We write
\begin{equation}\label{eqn3.1}
(\gamma)_*(\rho)
\end{equation}
to denote $\gamma\rho$ together with this orientation.  Clearly
$(\gamma)_*$ is functorial, and preserves the relative orientation
of $\rho,\tau$ whenever $\tau$ is a codimension-one face of~$\rho$.

Once and for all, fix orientations on the standard cells $\sigma_T$.
We can now put orientations on all the cells of~$W$.

\begin{definition}\label{defone}
Let $\sigma$ be a cell in $W_T$ with $\sigma = \gamma_0 \sigma_T$.
Let $O$ be the right $\bar\Gamma_T$-orbit in $\P^{n-1}$ corresponding
to~$\sigma$ as in Proposition~\ref{proptwo}.  If $O$ is orientable, we
give~$\sigma$ the orientation
\begin{equation}\label{eqn3.2}
(\text{orientation number of }\aaa) \cdot (\gamma_0)_*(\sigma_T).
\end{equation}
If $O$ is non-orientable, we give~$\sigma$ an arbitrary orientation.
\end{definition}

\begin{proposition}\label{propthree}
The quantity in~\eqref{eqn3.2} is well-defined.
\end{proposition}

\begin{proof}
 Let $\sigma$, $\gamma_0$, $O$, and $\aaa_0$ be as in
Definition \ref{defone}, with $O$ assumed orientable.  Assume $\sigma =
\gamma_1\sigma_T$ as well as $\gamma_0\sigma_T$.  Let $\aaa_1 =
\mathfrak{b}(\gamma_1)$.  Then $\gamma_1^{-1} \gamma_0 \in \bar\Gamma_T$.
By definition of $\bar\Gamma_T^+$, $(\gamma_0)_*(\sigma_T) =
(\gamma_1)_*(\sigma_T)$ if and only if $\gamma_1^{-1} \gamma_0 \in
\bar\Gamma_T^+$.  On the other hand, $\aaa_1 \gamma_1^{-1} \gamma_0
= [0:\dots:0:1] \gamma_0 = \aaa$, so $\aaa$ and $\aaa_1$ have
the same orientation number if and only if $\gamma_1^{-1} \gamma_0 \in
\bar\Gamma_T^+$.  \end{proof}

We must understand how $\bar\Gamma$ and $\Gamma$ act on the
orientations we have just chosen.  The following lemma is immediate
from \eqref{eqn3.1} and Definition~\ref{defone}.

\begin{lemma}\label{lemfour}
If $\sigma = \gamma_0 \sigma_T$ corresponds to an orientable~$O$, then
$\gamma_0$ carries the chosen orientation of $\sigma_T$ to the chosen
orientation of $\sigma$ times the orientation number of $\aaa$.
\end{lemma}

Here is a more general statement.

\begin{proposition}\label{propfour}
Let $\sigma = \gamma_0 \sigma_T$ and
$\sigma_1 = \gamma_1 \sigma$, for some $\gamma_0, \gamma_1\in
\bar\Gamma$.  Let $\aaa_0 = \mathfrak{b}(\gamma_0)$ and $\aaa_1 =
\mathfrak{b}(\gamma_1 \gamma_0)$.  Let $O_1, O$ be the right
$\bar\Gamma_T$-orbits in $\P^{n-1}$ containing $\aaa_1, \aaa_0$;
assume these orbits are both orientable.  Then $\gamma_1$ carries
$\sigma$ to $\sigma_1$ while multiplying the orientations by
\begin{equation}\label{eqn3.3}
(\text{orien.~number of $\aaa$})
     \cdot (\text{orien.~number of $\aaa_1$}).
\end{equation}
\end{proposition}

\begin{proof}
 Apply Lemma~\ref{lemfour} twice. 
\end{proof}

Fortunately, \eqref{eqn3.3} becomes trivial when we consider $\Gamma$
as opposed to $\bar\Gamma$.

\begin{proposition}\label{propfive} Let $\sigma_1 = \gamma_1 \sigma$
for $\gamma_1\in \Gamma$.  Let $\aaa_1, \aaa, O_1, O$ be as in
Proposition~\ref{propfour}, both orbits being assumed orientable.
Then $\gamma_1$ carries $\sigma$ to $\sigma_1$ with orientations
matching.
\end{proposition}

\begin{proof}
 We have $[0:\dots:0:1] \gamma_1 = [0:\dots:0:1]$ by the
definition of~$\Gamma$.  Hence $\aaa_1 = ([0:\dots:0:1] \gamma_1)
\gamma_0 = [0:\dots:0:1] \gamma_0 = \aaa$.  Thus the expression
in \eqref{eqn3.3} is a square either of~$+1$ or of~$-1$.  \end{proof}

\subsection{}\label{subsectfour} 

We compute the equivariant cohomology $H^*_\Gamma(W; \C)$ using a
spectral sequence, following the exposition in \cite[VII.7--8]{Br}.
(The spectral sequence there is for equivariant homology; we make the
appropriate modifications for cohomology.)

Let $\sigma$ be any cell in the well-rounded retract~$W$.  Recall that
$\Gamma_\sigma$ is the stabilizer of~$\sigma$ in~$\Gamma$.  Let
$\C_\sigma$ be the $\Gamma_\sigma$-module where $\gamma\in
\Gamma_\sigma$ acts by $+1$ if $\gamma$ preserves the orientation
of~$\sigma$ and by $-1$ if it does not.

For each~$i$, let $W_{(i)}$ be a fixed set of representatives of the
$\Gamma$-orbits of the cells in $W$ of dimension~$i$.  The $E_1$ term
of the spectral sequence is
\begin{equation}\label{eqn4.1}
E_1^{i,j} = \bigoplus_{o\in W_{(i)}} H^j(\Gamma_o; \C_o).
\end{equation}
These cells~$o$ (omicron) are in one-to-one correspondence with
the~$O$ of Proposition~\ref{proptwo}, as~$T$ runs through the types of
cells of dimension~$i$.  Because $\C$ is a field of characteristic
zero, all the terms in \eqref{eqn4.1} vanish when $j\ne 0$.  In
particular, the spectral sequence collapses at $E_2$.  The term
$H^0(\Gamma_o; \C_o)$ is the subset of $\Gamma_o$-invariants in the
module $\C_o$.

\begin{proposition}\label{propsix} 
For any cell $\sigma$, let~$T$ be the type of~$\sigma$, and let $O$ be
the right $\bar\Gamma_T$-orbit in $\P^{n-1}$ corresponding to~$\sigma$
as in Proposition~\ref{proptwo}.  Then $H^0(\Gamma_\sigma; \C)$
is~$\C$ if $O$ is orientable, and is~0 if $O$ is non-orientable.
\end{proposition}

\begin{proof}
 When $O$ is orientable, this follows from Proposition~\ref{propfive},
merely because $\Gamma_\sigma \subset \Gamma$.  Now assume $O$ is
non-orientable.  Let $\sigma = \gamma_0 \sigma_T$, with $\aaa =
\mathfrak{b}(\gamma_0)$.  As we have said above, there is some $\gamma_1
\in \bar\Gamma_T \smallsetminus \bar\Gamma_T^+$ with $\aaa \gamma_1 = \aaa$.
By~(2.2), $\Gamma_\sigma = (\gamma_0 \bar\Gamma_T \gamma_0^{-1}) \cap
\Gamma$.  The element $\gamma_0 \gamma_1 \gamma_0^{-1}$ is in $\Gamma$
by Lemma~2, so it is in $\Gamma_\sigma$; clearly it carries $\sigma$
to itself while reversing orientation.  Hence $\Gamma_\sigma$ acts
non-trivially on $\C_\sigma$, meaning $H^0(\Gamma_\sigma; \C_\sigma) =
0$.  \end{proof}

\begin{remark}
The proposition shows that the $\Gamma$-orbits of cells
coming from non-orientable~$O$ contribute nothing to our spectral
sequence.  We ignore these objects for the rest of the computation,
tacitly assuming that all $O$'s mentioned from now on are orientable.
\end{remark}

To summarize:

\begin{proposition}\label{propseven} 
The $E_1^{i,0}$ term of the equivariant
cohomology spectral sequence for $H^i_\Gamma(W; \C)$ is a direct sum
$\bigoplus_o \C$, where $o$ runs through a set of $i$-cells in
one-to-one correspondence with the orientable right
$\bar\Gamma_T$-orbits in $\P^{n-1}$, for all types~$T$ of cells of
dimension~$i$. 
\end{proposition}

\subsection{}\label{subsectfive}
 In \S\ref{subsectsix}, we will describe the boundary maps
$d_1$ of the spectral sequence.  These are the only differentials we
need consider, since the sequence collapses at~$E_2$.  In this
subsection, we give some details concerning how the cells meet at
their boundaries.

As usual, a \emph{facet} of a cell $\sigma$ is any face of $\sigma$
of codimension one.  Let $\mathcal{F}_\sigma$ be the set of facets
of~$\sigma$ in~$W$.

We will need to understand how $\mathcal{F}_\sigma$ breaks up into
orbits under the action of $\Gamma_\sigma$ (and to choose a set $\mathcal{F}_\sigma'$ of representatives of these orbits in
$\mathcal{F}_\sigma$).  We do this in Proposition~\ref{propeight}
below.  We will start by determining $\mathcal{F}_{\sigma_T}$ for the
standard cells $\sigma_T$ in a form suited to our computation.  We
will then determine $\mathcal{F}_\sigma$ for any~$\sigma$.

We make two conventions.  (i)~If types~$T$ and $T'$ occur in the same
discussion, it is assumed that a cell of type~$T$ has at least some
cells of type~$T'$ as facets.  (ii)~If $\sigma = x \sigma_T$ for some
$x\in \bar\Gamma$, we identify $\sigma$ with the coset $x\bar\Gamma_T$
as in Proposition~\ref{propone}.  In the expressions of the form
$$
\bigcup_{T'} (*),
$$
the ($*$) will be a finite union of cosets, say $x_1 \bar\Gamma_{T'}
\coprod \cdots \coprod x_k \bar\Gamma_{T'}$ for $x_1,\dots, x_k \in
\bar\Gamma$.  It is understood that~($*$) corresponds to the set of
cells $x_1\sigma_{T'}, \dots, x_k\sigma_{T'}$, as in
Proposition~\ref{propone}.
Even when~($*$) is a more complicated object, like a double coset, its
meaning is that one should decompose it into single cosets by choosing
appropriate representatives (which will be the $x_i$).

The boundary of $\sigma_T$ is a union of (the closures of) various
cells of type $T'$:
\begin{equation}\label{eqn5.1}
\mathcal{F}_{\sigma_T} =
\bigcup_{T'} \beta_{T',1} \bar\Gamma_{T'} \coprod \dots \coprod
\beta_{T',k} \bar\Gamma_{T'}
\end{equation}
for some $\beta_{T',i} \in \bar\Gamma$.  However, \eqref{eqn5.1} is invariant
under the left action of $\sigma_T$'s stabilizer $\bar\Gamma_T$.
Hence there must be finitely many $\alpha_{(T, T', \iota)} \in \bar\Gamma$
such that
\begin{equation}\label{eqn5.1a}
\mathcal{F}_{\sigma_T} =
\bigcup_{T'} \bar\Gamma_T \alpha_{(T, T', 1)}\bar\Gamma_{T'}
\coprod \dots \coprod
\bar\Gamma_T \alpha_{(T, T', k)} \bar\Gamma_{T'}.
\end{equation}

We have computed the $\alpha_{(T, T', \iota)}$ by hand for $n\le4$.
In our cases of interest ($n=4$, $\dim\sigma \ge4$), one finds there
is only one~$\iota$---that is, the right-hand side of~\eqref{eqn5.1a}
is actually just one double coset.  In fact, we find that we may take
$\alpha_{(T,T',\iota)}$ to be the identity except in one case, where
$(T,T') = ($5b, 4b$)$ in the notation of \cite{M91} and
$$
\alpha_{\text{5b}, \text{4b}, 1} =
\left (\begin{array}{cccc}
0&0&1&0 \\
0&0&0&1 \\
-1&0&1&0 \\
0&-1&1&0
\end{array}\right).
$$
From now on, we will write $\alpha$ for $\alpha_{(T, T', 1)}$, the
dependence on $T$ and $T'$ being understood.

To find an expression for $\mathcal{F}_\sigma$ for $\sigma = \gamma_0
\sigma_T$, we multiply~\eqref{eqn5.1a} by $\gamma_0$ to obtain
\begin{equation}\label{eqn5.2}
\mathcal{F}_\sigma = \bigcup_{T'}  \gamma_0 \bar\Gamma_T \alpha
\bar\Gamma_{T'}.
\end{equation}
To exhibit this as a union of cosets of the form $x_i\bar\Gamma_{T'}$,
we rewrite it as
\begin{align}
\mathcal{F}_\sigma
&= \bigcup_{T'}  \gamma_0 \alpha (\alpha^{-1} \bar\Gamma_T
                \alpha) \bar\Gamma_{T'} \label{eqn5.3a}\\
&= \bigcup_{T'}  \gamma_0 \alpha (\alpha^{-1} \bar\Gamma_T
                \alpha / (\alpha^{-1} \bar\Gamma_T
                \alpha \cap \bar\Gamma_{T'})) \bar\Gamma_{T'}\label{eqn5.3b},
\end{align}
The last formula (for each $T'$) is a disjoint union of single cosets,
in one-to-one correspondence with a set of representatives of the
cosets in $\alpha^{-1} \bar\Gamma_T \alpha / (\alpha^{-1} \bar\Gamma_T
\alpha \cap \bar\Gamma_{T'})$.  It is a matter of formal manipulation
to get an expression for $\Gamma_\sigma \bs \mathcal{F}_\sigma$:
\begin{align}
\Gamma_\sigma \bs \mathcal{F}_\sigma
&= \bigcup_{T'}  (\gamma_0 \bar\Gamma_T \gamma_0^{-1} \cap \Gamma) \bs
                \gamma_0 \alpha (\alpha^{-1} \bar\Gamma_T
                \alpha / (\alpha^{-1} \bar\Gamma_T
                \alpha \cap \bar\Gamma_{T'}))  \bar\Gamma_{T'} \\
&= \bigcup_{T'}  \gamma_0 (\bar\Gamma_T \cap \gamma_0^{-1} \Gamma
                                                      \gamma_0) \bs 
                \alpha (\alpha^{-1} \bar\Gamma_T
                \alpha / (\alpha^{-1} \bar\Gamma_T
                \alpha \cap \bar\Gamma_{T'}))  \bar\Gamma_{T'} \\
&= \bigcup_{T'} 
                \underbrace{
                \gamma_0 \alpha 
                [\underbrace{
                  \alpha^{-1} (\bar\Gamma_T \cap \gamma_0^{-1} \Gamma \gamma_0)
                        \alpha}_{A}
                \bs 
                \alpha^{-1} \bar\Gamma_T \alpha}_{B}
                 /
                \underbrace{(\alpha^{-1} \bar\Gamma_T
                \alpha \cap \bar\Gamma_{T'})}_{C}]
                 \bar\Gamma_{T'}\label{eqn5.3stuff}
\end{align}
where for each $T'$ the right-hand side of \eqref{eqn5.3stuff} is
again expressed as a disjoint union of $\bar\Gamma_{T'}$-cosets, in
one-to-one correspondence with a set of representatives of the double
coset expression in square brackets.

We must interpret \eqref{eqn5.3a}--\eqref{eqn5.3b} in terms of the
$\bar\Gamma_{T'}$-orbits in $\P^{n-1}$.  As usual, let $O$ be the
right $\bar\Gamma_T$-orbit corresponding to~$\sigma$, with $\aaa =
\mathfrak{b}(\gamma_0)$.  Let $\mathbf{c} = \mathfrak{b}(\gamma_0
\alpha)$, so that $\mathbf{c} = \aaa \alpha$.  Now the quantity $A$ in
\eqref{eqn5.3stuff} is exactly the subgroup of $\alpha^{-1}
\bar\Gamma_T \alpha$ that preserves $\mathbf{c}$.  Hence any set of
coset representatives for $B$ in \eqref{eqn5.3stuff} is a set of
matrices whose bottom rows, in $\P^{n-1}$, are the members of the
right $(\alpha^{-1} \bar\Gamma_T \alpha)$-orbit of $\mathbf{c}$.
Equivalently, any set of representatives for $B$ is a set of matrices
whose bottom rows are exactly the members of $O\cdot \alpha$.  The
group denoted $C$ in \eqref{eqn5.3stuff} acts on this orbit $O\cdot
\alpha$, decomposing it into suborbits; for given $T'$, the disjoint
$\bar\Gamma_{T'}$-cosets in~\eqref{eqn5.3a}--\eqref{eqn5.3b} are in
one-to-one correspondence with these suborbits.  We summarize this
result as follows:

\begin{proposition}\label{propeight}
 Let $T$, $T'$ and $\alpha$ be as introduced in this subsection.  Let
$\sigma$ be a cell of type~$T$, represented by the
$\bar\Gamma_T$-orbit $O$ in the manner of \S\ref{subsectthree}.  Decompose the
orbit $O\cdot \alpha$ into its suborbits $O_1,\dots,O_k$ under the
group $C = \alpha^{-1} \bar\Gamma_T \alpha \cap
\bar\Gamma_{T'}$.  Let $\aaa_j \in O_j$.  Let $\gamma_j \in
\bar\Gamma$ be chosen so that $\mathfrak{b}(\gamma_j) = \aaa_j$.
(We may, in fact, take $\gamma_j$ to be of the form $\gamma_0
\hat\gamma_j \alpha$ for some $\hat\gamma_j \in \Gamma_T$.)  Then the
union over all $T'$ of the cells
$$
\{ \gamma_1 \sigma_{T'}, \dots, \gamma_k \sigma_{T'} \}
$$
is a set $\mathcal{F}'_\sigma$ of representatives for $\Gamma_\sigma \bs \mathcal{F}_\sigma$.  
\end{proposition}

For all pairs $T, T'$, we store the intersection $C =
\alpha^{-1} \bar\Gamma_T \alpha \cap \bar\Gamma_{T'}$ in our program.

\subsection{}\label{subsectsix}
We can now determine the boundary maps $d_1$ in the
spectral sequence.  Recall that $W_{(i)}$ is a fixed set of
representatives of the $\Gamma$-orbits of the cells in $W$ of
dimension~$i$.  Superseding the use of~$o$ in \S\ref{subsectfour}, we let
$\sigma$ run through $W_{(i+1)}$, and let $\tau$ run through the set
$\mathcal{F}'_\sigma$ of representatives of the facets of~$\sigma$.  To
compute $H^i_\Gamma(W; \C)$ for subgroups of $\SL_4(\Z)$ for $i=5,6$,
we must work out the map $d_1^{i,0}$ for $i=4$ and~5.

We follow \cite[VII.8]{Br}, taking the dual to turn homology into
cohomology.  The map
\begin{equation}\label{thatline}
d_1^{i,0} : \bigoplus_{{\tau_0}\in W_{(i)}} H^0(\Gamma_{\tau_0}; \C_{\tau_0})
 \to \bigoplus_{\sigma\in W_{(i+1)}} H^0(\Gamma_\sigma; \C_\sigma)
\end{equation}
is a sum of terms, one for each pair $(\sigma,\tau)$; here $\tau_0$ is
the fixed representative in $W_{(i)}$ that is $\Gamma$-equivalent
to~$\tau$.  From now on, we focus on a single pair $\sigma$, $\tau$.
Call their types~$T$ and $T'$, respectively.  Let $O$ be the right
$\bar\Gamma_T$-orbit in $\P^{n-1}$ associated to~$\sigma$ by
Proposition~\ref{proptwo}.  If $\sigma = \gamma_0\sigma_T$, let $\aaa = \mathfrak{b}(\gamma_0) \in O$.  Let $\Gamma_{\sigma\tau} = \Gamma_\sigma \cap
\Gamma_\tau$.

As in \cite[p.~176]{Br}, the term for $(\sigma, \tau)$ in $d_1^{i,0}$
is the composition
\begin{equation}
\xymatrix{H^0(\Gamma_{\tau_0}; \C_{\tau_0})
  \ar[r]^{v_\tau}&
H^0(\Gamma_\tau; \C_\tau)
  \ar[r]^{u_{\sigma\tau}}&
H^0(\Gamma_{\sigma\tau}; \C_\sigma)
  \ar[r]^{t_{\sigma\tau}}&
H^0(\Gamma_\sigma; \C_\sigma)}.
\end{equation}

We now give the definition of these maps and show how to compute them.
Note that all the coefficient modules are copies of~$\C$ on which the
groups act trivially.

The map $t_{\sigma\tau}$ is the transfer map
$\C\to\C$ given by multiplication by the scalar $[\Gamma_\sigma :
\Gamma_{\sigma\tau}]$.  We evaluate this scalar as follows.  First of
all, $\Gamma_\sigma = (\gamma_0
\bar\Gamma_T \gamma_0^{-1}) \cap \Gamma$, which has the same
cardinality as $\bar\Gamma_T \cap (\gamma_0^{-1} \Gamma \gamma_0)$.
The latter is the subgroup of $\bar\Gamma_T$ that fixes $\aaa$.
Thus
\begin{equation}\label{eqn6.2}
\#(\Gamma_\sigma) = \frac
{\#(\bar\Gamma_T)}
{\#(\bar\Gamma_T\text{-orbit of } \aaa)}
= \frac{\#(\bar\Gamma_T)} {\#(O)}.
\end{equation}
The numerator is known because we stored $\bar\Gamma_T$.  The denominator
is easily recovered from the computer's lists of orbits.

We now evaluate the order of $\Gamma_{\sigma \tau}$.  As before,
$\Gamma_\sigma = \gamma_0 \bar\Gamma_T \gamma_0^{-1} \cap \Gamma$.  By
Proposition~\ref{propeight}, $\tau = \gamma_0 \hat\gamma \alpha
\sigma_{T'}$ for some $\hat\gamma \in \bar\Gamma_T$.  Hence
$$
\Gamma_\tau = ((\gamma_0 \hat\gamma \alpha) \bar\Gamma_{T'}
                 (\alpha^{-1} \hat\gamma^{-1} \gamma_0^{-1})) \cap \Gamma.
$$
Writing out $\Gamma_{\sigma\tau} = \Gamma_\sigma \cap \Gamma_\tau$ and
conjugating by $(\gamma_0 \hat\gamma \alpha)$, we find that
$\Gamma_{\sigma\tau}$ has the same cardinality as
\begin{equation}\label{eqn6.2a}
\underbrace{(\alpha^{-1} \bar\Gamma_T \alpha \cap
                               \bar\Gamma_{T'})}_{D}
  \cap
           (\alpha^{-1} \hat\gamma^{-1} \gamma_0^{-1}
             \Gamma \gamma_0 \hat\gamma \alpha)
\end{equation}
The group~$D$ is the same as the group $C$ from \eqref{eqn5.3stuff},
and is one of the groups we've stored.  And the whole group in
\eqref{eqn6.2a} is exactly the subgroup of $D$ that fixes the point
$\aaa_j$ of Proposition~\ref{propeight}.  So
\begin{equation}\label{eqn6.3}
\#(\Gamma_{\sigma\tau}) = \frac
{\#(\alpha^{-1} \bar\Gamma_T \alpha \cap \bar\Gamma_{T'})}
{\#(O_j)}.
\end{equation}
Again, the numerator is known from what's stored, and the denominator
is easy to evaluate.

We have
\begin{equation}\label{eqn6.4}
t_{\sigma\tau} = \frac{\text{\eqref{eqn6.2}}}{\text{\eqref{eqn6.3}}}.
\end{equation}
As a consistency check, the program signals an error if the
computed value of~\eqref{eqn6.4} is not an integer.

Brown's $u_{\sigma\tau} : \C \to \C$ is the composition
\[
H^0(\Gamma_\tau; \C_\tau) \to H^0(\Gamma_{\sigma\tau}; \C_\tau)
             \to H^0(\Gamma_{\sigma\tau}; \C_\sigma).
\]
The first arrow is induced from the inclusion $\Gamma_{\sigma\tau}
\hookrightarrow \Gamma_\tau$, and is easily seen to be the identity.
The second arrow is induced by the $\Gamma_{\sigma\tau}$-map
$\partial_{\sigma\tau} : \C_\sigma \to \C_\tau$, namely the
$(\sigma,\tau)$-component of the cellular boundary operator on~$W$.
Let $[\sigma:\tau]$ be $\pm1$ depending on whether the orientation on
$\sigma$ from Definition~\ref{defone} does or does not induce the
orientation on the facet $\tau$ from Definition~\ref{defone}.  Then
$\partial_{\sigma\tau}$ is the map $\C\to\C$ given by the scalar
$[\sigma:\tau]$.  Thus $u_{\sigma\tau} = [\sigma : \tau]$.

To evaluate $[\sigma:\tau]$, we reduce the problem to evaluating the
boundary operator on a small list of pairs of cells.  Write
$\tau=\gamma_0 \hat\gamma \alpha \sigma_{T'}$.  Then $\sigma =
\gamma_0 \hat\gamma \sigma_T$, since $\hat\gamma \in \bar\Gamma_T$.
If $\mathbf{x} \in \P^{n-1}$ is part of any oriented $\bar\Gamma_T$-orbit
(for any~$T$), write $\op{sgn}_T(\mathbf{x})$ for the orientation number
of $\mathbf{x}$ with respect to $\bar\Gamma_T$.  By Definition~\ref{defone},
$\sigma$ and $\tau$---together with their standard orientations---are
given as follows:
\begin{align}
\sigma &= (\gamma_0 \hat\gamma)_* (\sigma_T)
             \cdot \op{sgn}_T(\mathfrak{b}(\gamma_0 \hat\gamma)) \\
\tau &= (\gamma_0 \hat\gamma \alpha)_*(\sigma_{T'}) \cdot
\op{sgn}_{T'}(\mathfrak{b}(\gamma_0 \hat\gamma \alpha)).
\end{align}
Because $(\dots)_*$ is functorial,
$$
\tau = (\gamma_0 \hat\gamma)_* (\alpha)_* (\sigma_{T'})
             \cdot \op{sgn}_{T'}(\mathfrak{b}(\gamma_0 \hat\gamma
                                            \alpha)).
$$
Since $(\dots)_*$ preserves the $[\text{\ }:\text{\ }]$ relation, we
may cancel out $(\gamma_0\hat\gamma)_*$'s, obtaining
\begin{equation}\label{eqn6.5}
[\sigma:\tau] =
  \underbrace{\op{sgn}_T(\mathfrak{b}(\gamma_0 \hat\gamma))}_{E}
  \cdot
  \underbrace{\op{sgn}_{T'}(\mathfrak{b}(\gamma_0 \hat\gamma
                                       \alpha))}_{F}
  \cdot
  [\sigma_T : (\alpha)_*(\sigma_{T'})].
\end{equation}

We evaluate each factor in~\eqref{eqn6.5} in turn.  In~$F$,
$\mathfrak{b}(\gamma_0 \hat\gamma \alpha)$ is the point $\aaa_j$ of
Proposition~\ref{propeight}, and $\op{sgn}_{T'}(\aaa_j)$ is simply the
orientation number of $\aaa_j$ in its $\bar\Gamma_{T'}$-orbit.
Similarly, $E$ is the orientation number of $\mathfrak{b}(\gamma_0
\hat\gamma)$, which is a $\bar\Gamma_T$-translate of $\aaa$, in its
$\bar\Gamma_T$-orbit.  In practice, though, we do not know
$\hat\gamma$ or $\mathfrak{b}(\gamma_0 \hat\gamma)$ explicitly, so it
is easier to evaluate~$E$ by another method.  Note that
$\mathfrak{b}(\gamma_0 \hat\gamma) \cdot \alpha =
\mathfrak{b}(\gamma_0 \hat\gamma \alpha) = \aaa_j$, an element of the
$(\alpha^{-1} \bar\Gamma_T \alpha)$-orbit $O\cdot \alpha$ of
Proposition~\ref{propeight}.  It is easy to compute all the points in
$O\cdot \alpha$, using the group $\alpha^{-1} \bar\Gamma_T \alpha$
(which was stored).  For some $\mathbf{x}$, we will have the equation
$\aaa _{j} = \mathbf{x}\cdot \alpha $.  Then $E$ will be the
orientation number for this $\mathbf{x}$.

The quantities $[\sigma_T : (\alpha)_*(\sigma_{T'})]$ are evaluated by
hand, for all $T, T'$ occurring in~\eqref{eqn5.1a}.  Several issues
arise.  First, we must find $[\sigma_T : \sigma_{T'}]$ whenever the
standard cell $\sigma_{T'}$ is a facet of $\sigma_T$.  No matter how
we choose the orientations at the start, it is in general impossible
to arrange our choices so that all the relative orientations are
positive.  As a general illustration, if a 0-cell, two 1-cells, and a
2-cell meet locally in a picture like the first quadrant of $\R^2$,
and if both 1-cells are oriented to point away from the origin, then
the pair (2-cell, $x$-axis) must have relative orientation opposite to
that of the pair (2-cell, $y$-axis), no matter how we orient the
2-cell.  In the $\SL_4(\Z)$ case, one can draw a schematic picture of
how the standard cells meet and can read off all the relative
orientations $[\sigma_T : \sigma_{T'}]$.

Second, we must compute $(\alpha)_*(\sigma_{T'})$ when $\alpha$ isn't
the identity, knowing the orientation on $\sigma_{T'}$.  This involves
the same techniques as in Remark \ref{remone}.

Third, we must find some facet $\upsilon$ of $\sigma_T$ whose
orientation is known, and must compare $[\sigma_T : \upsilon]$ to
$[\sigma_T : (\alpha)_*(\sigma_{T'})]$.  In practice, if $\dim
\sigma_T = k$, this means finding a chain of $(k-1)$-cells between
$\upsilon$ and $(\alpha)_*(\sigma_{T'})$ such that consecutive members
of the chain meet in faces of dimension~$k-2$, and comparing the
orientations of $\upsilon$ and $(\alpha)_*(\sigma_{T'})$ across the
$(k-2)$-faces.

Finally, we must compute $v_\tau$, but this is easy.  It is induced by
the conjugation action of the element of $\Gamma$ that carries $\tau$
to $\tau_0$.  But by Proposition~\ref{propfive}, any element of
$\Gamma$ preserves the orientations on the cells.  Hence one finds
$v_\tau = +1$.

\section{Hecke operators}\label{heckeop.sect}
\subsection{}
We identify the cochain complex in \eqref{thatline} with the complex
of cellular cochains on $W$ by identifying an $i$-cell $\sigma $ with
a generator of $H^{0} (\Gamma _{\sigma }; \C_{\sigma })$, taking care
to make the signs match.  Formulas like $\sum n (\sigma )\sigma $ will
denote the corresponding cocycles in either complex, and will be
referred to as \emph{$W$-cocycles}  

Let $\beta \in H^{5} (\Gamma ; \C)$ be a class, and let $u=\sum n
(\sigma )\sigma$ be a representative for $\beta$ in terms of the
previous paragraph.  Let $T (l,k)$ be a Hecke operator.  To compute
the action of $T (l,k)$ on $u$, we do the following:
\begin{enumerate}
\item Convert $u$ to a reduced $1$-sharbly cycle $\xi $, and then
compute $T (l,k)(\xi)$ using \eqref{how.hecke.acts} in
\S\ref{that.section}\label{itemone}.
\item Use the algorithm from \cite{gunn} to write $T (l,k) (\xi )$
as a sum of reduced $1$-sharbly cycles\label{already}. 
\item Convert these reduced $1$-sharbly cycles to
$W$-cocycles\label{itemthree}.
\end{enumerate}

Step \ref{already} is described in detail in \cite{gunn}, and we refer
the reader to that article.  In this section, we focus on steps
\ref{itemone} and \ref{itemthree} in the context of
\S\ref{imp.details.sect}.

We begin with a definition from \cite{gunn}:
\begin{definition}\label{lifts}
\cite[Definition 5.3]{gunn} Let $\uu$ be a basis element of the
$k$-sharblies $S_{k}$.  Then a \emph{lift} for $\uu $ is an $n\times
(n+k)$ integral matrix $M$ with primitive columns such that
$[M_{1},\dots ,M_{n+k}] = \uu $, where $M_{i}$ is the $i$th column of
$M$.
\end{definition}

Modulo the action of $\GL_{4} (\Z )$, there is only one orbit of
reduced basis $0$-sharblies and only three orbits of reduced basis
$1$-sharblies.  The identity matrix serves as a lift for a member of
the first orbit, and lifts representing elements of the latter three
orbits are
\[
\left(\begin{array}{ccccc}
1&0&0&0&1\\
0&1&0&0&1\\
0&0&1&0&1\\
0&0&0&1&1
\end{array} \right),\quad 
\left(\begin{array}{ccccc}
1&0&0&0&1\\
0&1&0&0&1\\
0&0&1&0&1\\
0&0&0&1&0
\end{array} \right),\quad\hbox{and}\quad
\left(\begin{array}{ccccc}
1&0&0&0&1\\
0&1&0&0&1\\
0&0&1&0&0\\
0&0&0&1&0
\end{array} \right).
\]
We call these the \emph{standard} $0$- and $1$-sharblies.  The sets of
primitive vectors indexing the standard $6$ and $5$-cells in $W$ coincide
with the sets of column vectors of these matrices.  By abuse of
language we will speak of the ``standard sharbly of type $T$,'' and
will use the notation $\uu _{T}$.

\subsection{}\label{liftsect}
Given a sharbly cycle $\xi $, we denote by $\support \xi $
the support of $\xi $.  Suppose that $\Gamma $ is torsion-free.
Then according to \cite{gunn}, a $1$-sharbly cycle $\xi $ mod $\Gamma $
with coefficients in a ring $R$ can be encoded by a collection 
of $4$-tuples $(\uu , n (\uu ), \{\vv \}, \{L (\vv ) \})$,
where
\begin{enumerate}
\item $\uu \in \support \xi $,
\item $n (\uu )\in R$,
\item $\{\vv  \} = \support \partial \uu $, and  
\item $\{L (\vv ) \}$ is a $\Gamma $-equivariant set of lifts for $\{\vv \}$.  \label{stepfour}
\end{enumerate}
The $\Gamma $-equivariance condition in \ref{stepfour} is the
following.  Suppose that for $\uu ,\uu '\in \support \xi $ there exist
$\vv \in \support (\partial \uu ) $ and $\vv '\in \support (\partial
\uu ')$ such that $\vv = \gamma \cdot \vv '$ for some $\gamma \in
\Gamma $.  Then we require $L (\vv )= \gamma L (\vv ')$.

In the case under study, $\Gamma $ is \emph{not} torsion-free, and the
above data needs to be modified.  Suppose that a $0$-sharbly
$\vv \in \support \partial \uu $ has a nontrivial stabilizer $\Gamma
(\vv )\subset \Gamma $, and let $m$ be any lift of $\vv $.
Then in the cycle $\xi $ we replace $n (\uu )\uu $ with 
\[
\sum _{\gamma \in \Gamma (\vv )} \frac{n (\uu )}{\#\Gamma (\vv )}\uu
_{\gamma }, 
\]  
where $\uu _{\gamma }$ has the same data as $\uu $, except that we
give $\vv $ the lift $\gamma m$.  (Note that this is possible
in our case since the coefficient ring $R=\C$ is divisible.)

\subsection{}
Now we describe how to construct the data in \S\ref{liftsect} to
produce a $1$-sharbly chain $\xi $ corresponding to the $W$-cocycle
$u$.  There are two steps.

First, choose $\sigma$ such that $n (\sigma )\not =0$ in $u$.
According to Proposition \ref{propone}, the $5$-cell $\sigma \in
\Gamma \bs W$ is encoded as a coset $[\gamma \Gamma _{T}]$, where
$\gamma \in \Gamma $ and $\Gamma _{T}$ is the stabilizer of the
standard cell of type $T$.  Moreover, the coset $[\gamma \Gamma _{T}]$
is encoded as the orbit $O$, which can be regarded a set of triples
$\{(\aa, \pm 1, \gamma_{\aa })\}$, where $\aa \in \Proj ^{3} (\Z
/{N}\Z )$, $\gamma_{\aa}\in \Gamma _{T}$, and $\pm 1$ is the
orientation number (\S\ref{subsectthree}).  From $O$ we arbitrarily
choose a triple with orientation number $1$, and then using Hermite
normal form construct a matrix $\gamma \in \SL_{4} (\Z )$ with bottom
row equal to $\aa$.  Then the contribution of $\sigma $ to the
$1$-sharbly chain is
\[
n (\sigma ) \gamma \uu _{T}.
\]
We do this for all of $\support u$ and sum to produce $\xi $.  For
each $\uu \in \support \xi $, we write $\uu (\gamma )$ if we want to
indicate the element $\gamma $ used in the construction of $\uu $.

\subsection{}
At this stage, we have a $1$-sharbly chain, and we need to choose
lifts to reflect the cycle structure of $\xi $.  This we do as
follows.  In the spirit of \S\ref{subsectfive}, for each type $T$ we
choose a set of matrices $\Omega _{T}\subset \SL_{4} (\Z )$ such that
$\partial \colon S_{1}\rightarrow S_{0}$ can be written as 
\[
\partial \colon \uu _{T} \longmapsto \sum_{\omega \in \Omega
_{T}}\omega \vv ,
\]
where $\vv $ is the standard basis $0$-sharbly.  Note the absence of
signs in this map---the signs in the boundary map in Definition
\ref{sharbly.complex} have been encoded in the $\omega $'s, which may
nontrivially permute the column vectors of $\vv $.

Form the $0$-sharbly chain
\begin{equation}\label{bdy}
\sum _{\uu (\gamma ),\omega \in \Omega _{T}} \gamma \omega \vv ,
\end{equation}
where $\vv $ is the standard $0$-sharbly, we sum over all $\uu (\gamma
)\in \support \xi $, and $T$ is the type of $\uu (\gamma )$.  (In
\eqref{bdy} we have abbreviated $\alpha _{(T,T',i)}$ to $\alpha _{i}$,
since $T$ is determined by $\uu (\gamma )$, and there is only one type
of reduced $0$-sharbly mod $\SL _{4} (\Z )$.)  Note that this sum is
in $S_{*}$, not $\coinv{*}$; the only relations we apply are those
in the sharbly complex.

\subsection{}
After summing, we find that some $0$-sharblies cancel, and some
remain.  For those that canceled, we can choose any lifts we like, as
long as we choose the same lifts for all terms that cancel each other.

The remaining $0$-sharblies form a chain $\eta $ that vanishes in
$\coinv{0}\otimes \C$, and we must choose nontrivial lifts for them.
To do this, first arbitrarily choose lifts for each $0$-sharbly in
$\support \eta $.  The data we computed in \S\ref{subsecttwo} allows us to
easily compute the distinct orbits of $\Gamma $ in $\support \eta $.
We do this and order each orbit.  

Suppose $\vv _{0}$ is
the first $0$-sharbly in one of these orbits, and that it corresponds
to the triple $(\aaa _{0}, \pm 1, \gamma _{\aaa_{0} })$.  Let $\vv
_{0} (\gamma )$ be its lift.  Then if $\vv $ is any other $0$-sharbly
in $\vv_{0} $'s orbit, corresponding to the triple $(\aaa, \pm 1,
\gamma _{\aaa })$, we replace the lift $\vv (\gamma )$ by
\[
\vv (\gamma ) \longleftarrow \vv (\gamma ) \gamma _{\aaa }^{-1} \gamma
_{\aaa _{0}}
\]

After these lifts are constructed, the cycle $\xi $ is ready for input
in the Hecke operator program.

\subsection{}
Upon completion, the Hecke operator program returns a reduced
$1$-sharbly cycle, which we must convert to a $W$-cocycle.
So let $\xi = \sum
n (\uu )\uu $ be a reduced $1$-sharbly cycle, and let $\uu \in \support
\xi $.  First we determine which of the three types $T$ of standard
reduced $1$-sharblies $\uu $ has.  Then we must find a matrix $\gamma
(\uu )\in \SL _{4} (\Z )$ such that
\[
\gamma (\uu )\cdot\uu _{T} = \uu , 
\]
where $\uu _{T}$ is the standard $1$-sharbly of type $T$.  This is
straightforward, although one must be careful to incorporate the
orientation number of $\gamma (\uu )$.

In practice, the main step is the following.  Let $\uu $ be a reduced
$1$-sharbly basis element with lift $M$.  We choose a nonsingular
$4\times 4$ minor $m$ of $M$ and construct $m^{-1}$.  Then
$m^{-1}\cdot \uu $ will be a $1$-sharbly with lift $m^{-1}M$, and will
be standard except possibly for one column vector.  By multiplying
$m^{-1}M$ on the left by elements of the stabilizer of the standard
$0$-sharbly, we can eventually produce the standard $1$-sharbly with
same type as $\uu $.  This allows us to construct $\gamma (\uu )$.

\section{Numerical
results}\label{num.results.sect}
\subsection{}
In this section we present numerical data from our experiments.  As
mentioned in \S\ref{introduction}, to avoid floating point problems
with $\C$-coefficients we usually work with $\F = \F_{31991}$, the
finite field with $31991$ elements, and in some cases work with $\Z $
or $\Q $.  The computations were carried out on a variety of Unix
machines at Columbia and Oklahoma State.  The code for the cohomology
of $W$ (\S\ref{imp.details.sect}) was written in Common Lisp.  The
Hecke operator code (\S\ref{heckeop.sect}) was written in C++, and
used the LiDIA library \cite{lidia}.  Perl scripts patched together
the outputs of the various programs and produced the tables in
\S\ref{hecke.table.sect}.

\subsection{}\label{5.2sect}
We first describe how we performed linear algebra on the large sparse
matrices that arise in our computations.  Fix the level~$N$.  We use
the notation of \S\ref{subsectsix}, working over $\C$ at first, and
letting $V_{(i)}$ be the domain of the map $d^{i,0}_1$ in
\eqref{thatline}.  To compute $H^5(\Gamma; \C)$, we must find the
kernel of $d^{5,0}_1$ modulo the image of $d^{4,0}_1$.  We prefer to
find the kernel of a single matrix.  Regard the $d^{i,0}_1$ as
matrices acting on the left on column vectors, and let $\mathfrak{D}_{\Z }$
be the matrix where $d^{5,0}_1$ is stacked on top of the transpose of
$d^{4,0}_1$:
$$
\mathfrak{D}_{\Z } = \left (\begin{array}{c}
     d^{5,0}_1
 \\
     (d^{4,0}_1)^{tr}
\end{array}\right)
$$             
This matrix defines a map $V_{(5)} \to V_{(6)} \oplus V_{(4)}$, where
we use the standard inner product to identify the transpose of
$d^{4,0}_1$ with its adjoint.  The kernel of this map is the space of
\emph{harmonic 5-cocycles}; it is isomorphic to $H^5(\Gamma; \C)$.
Let $s = \dim V_{(5)}$, the number of columns of $\mathfrak{D}_{\Z }$.

The matrix $\mathfrak{D}_{\Z }$ has coefficients in $\Z$.  For any
ring $R$, set $\mathfrak{D}_{R} = \mathfrak{D}_{\Z }\otimes R$, and set
$\mathfrak{D} = \mathfrak{D}_{\F }$.  In the tables in
\S\ref{table.sect}, the value of ``rank'' we report in the rows
labeled $R$ is a number almost certainly equal to $\dim \ker
\mathfrak{D}_{R}$, whose computation is explained below.  This number is
also almost certainly equal to $\dim H^5(\Gamma; \C)$.

\subsection{}\label{lancz.sect}
When $\mathfrak{D}$ is very large, we could not have found its kernel
without a sparse version of the Lanczos algorithm.  This algorithm is
usually used with real or complex matrices, particularly for
eigenvalue problems.  Following ideas in \cite{ody}, we translated it
into the mod-$p$ setting.  Let $\mathfrak{E} = \mathfrak{D}^{tr} \cdot
\mathfrak{D}$, a symmetrized version of $\mathfrak{D}$.  We choose a
random non-zero seed vector $\vv$ with coefficients in $\F$ and
consider the sequence $\vv, \mathfrak{E}\vv, \mathfrak{E}^2 \vv,
\dots$.  The Lanczos algorithm shows us how to compute not this
sequence, but the sequence $\vv = \mathbf{q}_0, \mathbf{q}_1, \dots$
resulting from it by the Gram-Schmidt orthogonalization process.  We
perform the Gram-Schmidt process mod~$p$ in the naive way, using $\sum
x_j y_j$ for the inner product.  This means we don't have the usual
guarantee that $\sum x_j^2$ will be non-zero when $(x_1, x_2, \dots)
\ne (0, 0, \dots)$.  If the inner product is ever 0 for non-trivial
$(x_j)$, we simply abort and choose another random seed $\vv$; even
for large $\mathfrak{D}$, these aborts happen less than half the
time.

The strength of the algorithm is that the RAM only has to hold the
sparse matrix $\mathfrak{D}$ and a few vectors of storage.  It does
not have to hold $\mathfrak{E}$, which is dense in general.  The
$\mathbf{q}$'s form a dense matrix, but they may be stored on the
disk, not in RAM.  (In our implementation, $\mathfrak{D}^{tr}$ was
stored in RAM along with $\mathfrak{D}$.)

Let $k$ be the largest value for which the set $\{\qq_0, \dots,
\qq_{s-k}\}$ is linearly dependent.  Reading the $\qq$'s back in from
disk, the algorithm allows us to backsolve for a non-zero vector $\yy
\in \ker \mathfrak{E}$.  What we want is an element of $\ker
\mathfrak{D}$.  A priori, we only know $\ker \mathfrak{E} \supseteq
\ker \mathfrak{D}$, and we will see that the containment is not always
an equality (though it would be over $\R$ or $\C$).  However, by
checking $\mathfrak{D} \cdot \yy = 0$ directly, we always find in
practice that the $\yy$ we compute lie in $\ker \mathfrak{D}$.

Thus each successful run of our algorithm produces one kernel vector
for $\mathfrak{D}$.  It also produces~$k$.  One can easily show $k\geq
\dim\ker\mathfrak{E}$.  We run Lanczos up to 30 or 40 times with
different random seeds, and we find $k$ is independent of the random
seed used.\footnote{For one level, random seeds produced a certain
value of $k$, while one random seed produced a value that was greater
by $1$.  The Lanczos method behaves this way when $E =
\{\mathfrak{E}\vv, \mathfrak{E}^{2}\vv, \mathfrak{E}^{3}\vv,\dots \}$
does not span the image of $\mathfrak{E}$, but only a proper subspace
of the image.  Since $\vv$ is chosen randomly, it is extremely rare
for $E$ to span less than the full image; our data bears this
statement out.} We conclude that $\dim \ker \mathfrak{E} = k$; though
we have not proved this, the computational evidence seems conclusive.
It is clear that $\dim\ker\mathfrak{E}\geq \dim\ker\mathfrak{D}$.

To find a basis of $\ker \mathfrak{D}$, we run Lanczos many times
until we have a set $S$ of $k + 10$ elements of $\ker \mathfrak{D}$.
We use mod~$p$ Gram-Schmidt on subsets $S' \subseteq S$ to find
maximal linearly independent subsets of~$S$.  We start with several
different $S'$'s.  In all the cases we checked, we found that maximal
linearly indepdent sets in $S$ had a common cardinality~$k'$, that any
$k'$-element subset of $S$ was linearly independent, and that any
subset of $S$ with more than $k'$ elements was dependent.  We conclude
that $k'$ is the value of $\dim \ker \mathfrak{D}$; again, the
computational evidence is convincing, though all we have proved is $k'
\le \dim \ker \mathfrak{D}$.  In the rows marked "Lanczos" in
Table~~\ref{betti}, what we report as ``rank'' is~$k'$.  In the rows
marked $\Z $ or $\Q $, what we report (namely $k'$) is provably the
rank of $H^5(\Gamma , \C)$.  We do find $k' < k$ sometimes in
practice.

\subsection{}\label{table.sect} 
In Table~\ref{betti}, we give the results of our Betti number
computations.  {\em Gau\ss} means that ordinary Gaussian elimination
was used to find the kernel of $\mathfrak{D}$, and {\em Lanczos} means
the algorithm of \S\ref{lancz.sect} was used.  The entries marked with
$*$ are those for which $k=k'$.

\begin{table}[ht]
\begin{center}
\begin{tabular}{|c|c|c|r||c|c|c|r|}
\hline
Level&Coefficients&rank&remarks&Level&Coefficients&rank&remarks\\
\hline\hline
11&$\Z  $ & 2 &Gau\ss&33&$\F$ & 10* &Lanczos\\
12&$\F$ & 0 &"&34&$\F$ & 12* &"\\
13&$\Z  $ & 1 &"&35&$\F$ & 7*  &"\\
14&$\F$ & 2 &"&36&$\F$ & 24 &"\\
15&$\F$ & 2 &"&37&$\F$ & 8  &Gau\ss\\
16&$\F$ & 3 &"&38&$\F$ & 14* &Lanczos\\
17&$\Z  $ & 3 &"&39&$\F$ & 10* &"\\
18&$\F$ & 9 &"&40&$\F$ & 9* &"\\
19&$\Z  $ & 3 &"&41&$\F$ & 9*  &"\\
20&$\F$ & 2 &"&42&$\F$ & 17 &"\\
21&$\F$ & 3 &"&43&$\F$ & 10* &"\\
22&$\F$ & 7 &"&44&$\F$ & 18 &"\\
23&$\Z  $ & 5 &"&45&$\F$ & 27 &"\\
24&$\F$ & 2 &"&46&$\F$ & 19* &"\\
25&$\F$ & 7 &"&47&$\F$ & 11* &"\\
26&$\F$ & 7 &"&48&$\F$ & 26 &"\\
27&$\F$ & 12 &"&49&$\F$ & 33* &"\\
28&$\F$ & 7 &"&50&$\F$ & 34 &"\\
29&$\F, \Q $ & 6 &"&51&$\F$ & 19*&"\\
30&$\F$ & 8 &Lanczos&52&$\F$ & 21* &"\\
31&$\F$ & 6  &Gau\ss&53&$\F$ & 17* &"\\
32&$\F$ & 12 &Lanczos&&&&\\
\hline
\end{tabular}
\end{center}
\caption{Probable Betti numbers for $H^{5} (\Gamma _{0} (N); \C )$.}\label{betti}
\end{table}

\subsection{}\label{hecke.table.sect}
Next, we present the Hecke data we computed.  These computations are
much more arduous than computing Betti numbers, and grow in complexity
very fast as a function of the number of cells of the retract $W$ mod
$\Gamma $.  Hence we were able to compute only a few Hecke operators,
usually only $T (2,*)$ and $T (3,*)$.  Beyond level $20$, it becomes
infeasible to compute $T (3, *)$; hence most of our data at large
levels is only for $T (2,*)$.  Happily this is usually sufficient to
guess persuasively what is happening with the cohomology.

\newcommand{\tott}{${\text{IIa}}$}
\newcommand{\ttoo}{${\text{IIb}}$}
\newcommand{\foot}{${\text{IV}}$}
\newcommand{\testentry}{\medskip\hrule\smallskip}

We use the following conventions and abbreviations.  All polynomials
should be considered as elements of $\bar\F[X]$, where the bar denotes
algebraic closure.  We denote the $p$-adic cyclotomic character of
$G_\Q$ by $\epsilon$, so that $\epsilon(\Frob_l) = l$ for any $l \ne
p$.  The symbol \tott\ (resp. \ttoo, \foot) denotes a Galois
representation of the form $\epsilon^a \sigma_k \oplus \epsilon^b
\oplus \epsilon^c$, where $(k,a,b,c)$ is $(2,0,2,3)$
(resp. $(2,2,0,1)$, $(4,0,1,2)$), if $\sigma _{k}$ is the Galois
representation associated to a weight $k$ classical holomorphic
cuspidal newform $f\in S_{k}^{\text{new}} (N')$, where $N'$ divides
$N$.  The same symbol prefixed with an "E" denotes a Galois representation of
the same form except that $\sigma_k$ is the Galois representation attached
to an Eisenstein series of weight $k$ for $\Gamma_1(N')$.

The individual tables are organized as follows.  For each level we
give the rank, as defined at the end of \S\ref{lancz.sect}.  Then each
block gives data for the Hecke eigenspaces.  The first column gives
the type Galois representation seemingly attached to this eigenspace.
The second column gives the dimension of this eigenspace.  The third
column gives the index $l$ of the Hecke operator, and the fourth
column gives the corresponding factored \emph{Hecke polynomial}.  This
is the polynomial defined on the right of \eqref{star} in Definition
\ref{attached.def}; it succinctly encodes the Hecke action on any
vector in the eigenspace.  At the bottom of each table, we indicate
the $\Q $-splitting of the spaces of newforms $S_{k}=
S_{k}^{\text{new}} (N)$ under the action of the Hecke operators; this
data is from \cite{stein}.  After each table, we comment of the
eigenclasses.

\testentry
\newcommand{\rrrank}{\text{rank}}
%
%
%
%
\begin{center}
\begin{tabular}{|p{40pt}|p{40pt}|p{40pt}|p{200pt}|}
\hline
\multicolumn{4}{|l|}{Level \textbf{11}. $\rrrank = 2$.}\\
\hline
\hline
\tott&1&$T_2$&$ (1-4 X)(1-8 X)(1 +2 X + 2 X^2)$\\
&&$T_3$&$ (1-9 X)(1-27 X)(1 + X + 3 X^2) $\\
&&$T_5$&$ (1-25 X)(1-125 X)(1 - X + 5 X^2) $\\
&&$T_7$&$ (1-49 X)(1-343 X)(1 +2 X + 7 X^2) $\\
\hline
\ttoo&1&$T_2$&$ (1- X)(1-2 X)(1 +8 X + 32 X^2)$\\
&&$T_3$&$ (1- X)(1-3 X)(1 +9 X + 243 X^2) $\\
&&$T_5$&$ (1- X)(1-5 X)(1 -25 X + 3125 X^2) $\\
&&$T_7$&$ (1- X)(1-7 X)(1 +98 X + 16807 X^2) $\\
\hline
\multicolumn{4}{|l|}{$\dim S_2(11) = 1$, $ \dim S_4(11) = 2$}\\
\hline
\end{tabular}
\end{center}


{The weight 4 newform doesn't lift.}

\testentry
%
%
%
%
\begin{center}
\begin{tabular}{|p{40pt}|p{40pt}|p{40pt}|p{200pt}|}
\hline
\multicolumn{4}{|l|}{Level \textbf{13}. $\rrrank = 1$.}\\
\hline
\hline
\foot&1&$T_2$&$ (1- 2 X)(1-4 X)(1 +5 X + 8 X^2)$\\
&&$T_3$&$ (1- 3 X)(1-9 X)(1 +7 X + 27 X^2) $\\
&&$T_5$&$ (1- 5 X)(1-25 X)(1 +7 X + 125 X^2) $\\
&&$T_7$&$ (1- 7 X)(1-49 X)(1 +13 X + 343 X^2) $\\
\hline
\multicolumn{4}{|l|}{$\dim S_2(13) = 0$, $ \dim S_4(13) = 1+2$}\\
\hline
\end{tabular}
\end{center}


{Only the rational weight 4 newform lifts.}

\testentry
%
%
%
%
\begin{center}
\begin{tabular}{|p{40pt}|p{40pt}|p{40pt}|p{200pt}|}
\hline
\multicolumn{4}{|l|}{Level \textbf{14}. $\rrrank = 2$.}\\
\hline
\hline
\tott&1&$T_3$&$ (1-9 X)(1-27 X)(1 +2 X + 3 X^2)$\\
\hline
\ttoo&1&$T_3$&$ (1- X)(1-3 X)(1 +18 X + 243 X^2)$\\
\hline
\multicolumn{4}{|l|}{$\dim S_2(14) = 1$, $ \dim S_4(14) = 1+1$}\\
\hline
\end{tabular}
\end{center}


{No weight four newforms lift.}

\testentry
%
%
%
%
\begin{center}
\begin{tabular}{|p{40pt}|p{40pt}|p{40pt}|p{200pt}|}
\hline
\multicolumn{4}{|l|}{Level \textbf{15}. $\rrrank = 2$.}\\
\hline
\hline
\tott&1&$T_2$&$ (1-4 X)(1-8 X)(1 + X + 2 X^2)$\\
\hline
\ttoo&1&$T_2$&$ (1- X)(1-2 X)(1 +4 X + 32 X^2)$\\
\hline
\multicolumn{4}{|l|}{$\dim S_2(15) = 1$, $ \dim S_4(15) = 1+1$}\\
\hline
\end{tabular}
\end{center}


{No weight four newforms lift.}

\testentry
%
%
%
%
\begin{center}
\begin{tabular}{|p{40pt}|p{40pt}|p{40pt}|p{200pt}|}
\hline
\multicolumn{4}{|l|}{Level \textbf{17}. $\rrrank = 2$.}\\
\hline
\hline
\tott&1&$T_2$&$ (1-4 X)(1-8 X)(1 + X + 2 X^2)$\\
&&$T_3$&$ (1-9 X)(1-27 X)(1  + 3 X^2) $\\
&&$T_5$&$ (1-25 X)(1-125 X)(1 +2 X + 5 X^2) $\\
&&$T_7$&$ (1-49 X)(1-343 X)(1 -4 X + 7 X^2) $\\
\hline
\ttoo&1&$T_2$&$ (1- X)(1-2 X)(1 +4 X + 32 X^2)$\\
&&$T_3$&$ (1- X)(1-3 X)(1  + 243 X^2) $\\
&&$T_5$&$ (1- X)(1-5 X)(1 +50 X + 3125 X^2) $\\
&&$T_7$&$ (1- X)(1-7 X)(1 -196 X + 16807 X^2) $\\
\hline
\multicolumn{4}{|l|}{$\dim S_2(17) = 1$, $ \dim S_4(17) = 1+3$}\\
\hline
\end{tabular}
\end{center}


{No weight four newforms lift.}

\testentry
%
%
%
%
\begin{center}
\begin{tabular}{|p{40pt}|p{40pt}|p{40pt}|p{200pt}|}
\hline
\multicolumn{4}{|l|}{Level \textbf{19}. $\rrrank = 3$.}\\
\hline
\hline
\tott&1&$T_2$&$ (1-4 X)(1-8 X)(1  + 2 X^2)$\\
&&$T_3$&$ (1-9 X)(1-27 X)(1 +2 X + 3 X^2) $\\
&&$T_5$&$ (1-25 X)(1-125 X)(1 -3 X + 5 X^2) $\\
\hline
\ttoo&1&$T_2$&$ (1- X)(1-2 X)(1  + 32 X^2)$\\
&&$T_3$&$ (1- X)(1-3 X)(1 +18 X + 243 X^2) $\\
&&$T_5$&$ (1- X)(1-5 X)(1 -75 X + 3125 X^2) $\\
\hline
\foot&1&$T_2$&$ (1- 2 X)(1-4 X)(1 +3 X + 8 X^2)$\\
&&$T_3$&$ (1- 3 X)(1-9 X)(1 +5 X + 27 X^2) $\\
&&$T_5$&$ (1- 5 X)(1-25 X)(1 +12 X + 125 X^2) $\\
\hline
\multicolumn{4}{|l|}{$\dim S_2(19) = 1$, $ \dim S_4(19) = 1+3$}\\
\hline
\end{tabular}
\end{center}


{The rational weight four newform lifts.}

\testentry
%
%
%
%
\begin{center}
\begin{tabular}{|p{40pt}|p{40pt}|p{40pt}|p{200pt}|}
\hline
\multicolumn{4}{|l|}{Level \textbf{20}. $\rrrank = 2$.}\\
\hline
\hline
\tott&1&$T_3$&$ (1-9 X)(1-27 X)(1 +2 X + 3 X^2)$\\
\hline
\ttoo&1&$T_3$&$ (1- X)(1-3 X)(1 +18 X + 243 X^2)$\\
\hline
\multicolumn{4}{|l|}{$\dim S_2(20) = 1$, $ \dim S_4(20) = 1$}\\
\hline
\end{tabular}
\end{center}


{The rational weight four newform doesn't lift.}

\testentry
%
%
%
%
\begin{center}
\begin{tabular}{|p{40pt}|p{40pt}|p{40pt}|p{200pt}|}
\hline
\multicolumn{4}{|l|}{Level \textbf{21}. $\rrrank = 3$.}\\
\hline
\hline
\tott&1&$T_2$&$ (1-4 X)(1-8 X)(1 + X + 2 X^2)$\\
\hline
\ttoo&1&$T_2$&$ (1- X)(1-2 X)(1 +4 X + 32 X^2)$\\
\hline
\foot&1&$T_2$&$ (1- 2 X)(1-4 X)(1 +3 X + 8 X^2)$\\
\hline
\multicolumn{4}{|l|}{$\dim S_2(21) = 1$, $ \dim S_4(21) = 1+1+2$}\\
\hline
\end{tabular}
\end{center}


{Of the weight four rational newforms, only one lifts.}

\testentry
%
%
%
%
\begin{center}
\begin{tabular}{|p{40pt}|p{40pt}|p{40pt}|p{200pt}|}
\hline
\multicolumn{4}{|l|}{Level \textbf{23}. $\rrrank = 5$.}\\
\hline
\hline
\tott&2&$T_2$&$ (1-4 X)(1-8 X)(1 -\alpha X + 2 X^2)$\\
\hline
\ttoo&2&$T_2$&$ (1- X)(1-2 X)(1 -4 \alpha X + 32 X^2)$\\
\hline
\foot&1&$T_2$&$ (1- 2 X)(1-4 X)(1 +2 X + 8 X^2)$\\
\hline
\multicolumn{4}{|l|}{$\dim S_2(23) = 2$, $ \dim S_4(23) = 1+4$}\\
\hline
\end{tabular}
\end{center}


{Here $\alpha$ satisfies $\alpha^2-\alpha+1 = 0$.  The rational weight four newform lifts.}

\testentry
%
%
%
%
\begin{center}
\begin{tabular}{|p{40pt}|p{40pt}|p{40pt}|p{200pt}|}
\hline
\multicolumn{4}{|l|}{Level \textbf{25}. $\rrrank = 7$.}\\
\hline
\hline
\foot&1&$T_2$&$ (1- 2 X)(1-4 X)(1 + X + 8 X^2)$\\
\hline
EIIa&2&$T_2$&$(1-4X) (1-8X) (1+\alpha X + 2X^{2})$\\
\hline
EIV&2&$T_2$&$(1-2X) (1-4X) (1+\beta X + 8X^{2})$\\
\hline
EIIb&2&$T_2$&$(1-X) (1-2X) (1+\gamma X + 32X^{2})$\\
\hline
\multicolumn{4}{|l|}{$\dim S_2(25) = 0$, $ \dim S_4(25) = 1+1+1$}\\
\hline
\end{tabular}
\end{center}


{Of the three weight four rational newforms, only one lifts.  Here
$\alpha $ satisfies $\alpha ^{2}+1=0$, $\beta$ satisfies $\beta
^{2}+49=0$, and $\gamma $ satisfies $\gamma ^{2} + 16 = 0$.}

\testentry
%
%
%
%
\begin{center}
\begin{tabular}{|p{40pt}|p{40pt}|p{40pt}|p{200pt}|}
\hline
\multicolumn{4}{|l|}{Level \textbf{27}. $\rrrank = 12$.}\\
\hline
\hline
\tott&1&$T_2$&$ (1-4 X)(1-8 X)(1  + 2 X^2)$\\
\hline
\ttoo&1&$T_2$&$ (1- X)(1-2 X)(1  + 32 X^2)$\\
\hline
\foot&1&$T_2$&$ (1- 2 X)(1-4 X)(1 +3 X + 8 X^2)$\\
\hline
EIIa&3&$T_2$&$(1-4X) (1-8X)(1+X) (1+2X) $\\
\hline
EIV&3&$T_2$&$(1-2X)(1-4X)(1+X) (1+8X)$\\
\hline
EIIb&3&$T_2$&$(1-X)(1-2X)(1+4X) (1+8X) $\\
\hline
\multicolumn{4}{|l|}{$\dim S_2(27) = 1$, $ \dim S_4(27) = 1+1+2$}\\
\hline
\end{tabular}
\end{center}


{Of the two weight four rational newforms, only one lifts.}

\testentry
%
%
%
%
\begin{center}
\begin{tabular}{|p{40pt}|p{40pt}|p{40pt}|p{200pt}|}
\hline
\multicolumn{4}{|l|}{Level \textbf{29}. $\rrrank = 6$.}\\
\hline
\hline
\tott&2&$T_2$&$ (1-4 X)(1-8 X)(1 -\alpha X + 2 X^2)$\\
&&$T_3$&$ (1-9 X)(1-27 X)(1 +\alpha X + 3 X^2) $\\
\hline
\ttoo&2&$T_2$&$ (1- X)(1-2 X)(1 -4 \alpha X + 32 X^2)$\\
&&$T_3$&$ (1- X)(1-3 X)(1 +9 \alpha X + 243 X^2) $\\
\hline
\foot&2&$T_2$&$ (1- 2 X)(1-4 X)(1 -\alpha X + 8 X^2)$\\
&&$T_3$&$ (1- 3 X)(1-9 X)(1 - (3\alpha +8) X + 27 X^2) $\\
\hline
\multicolumn{4}{|l|}{$\dim S_2(29) = 2$, $ \dim S_4(29) = 2+5$}\\
\hline
\end{tabular}
\end{center}


{Here $\alpha$ satisfies $\alpha^2+2\alpha-1=0$.  In this example, the weight two and weight four newforms that lift are defined over the same quadratic extension of $\Q$.}

\testentry
%
%
%
%
\begin{center}
\begin{tabular}{|p{40pt}|p{40pt}|p{40pt}|p{200pt}|}
\hline
\multicolumn{4}{|l|}{Level \textbf{31}. $\rrrank = 6$.}\\
\hline
\hline
\tott&2&$T_2$&$ (1-4 X)(1-8 X)(1 -\alpha X + 2 X^2)$\\
&&$T_3$&$ (1-9 X)(1-27 X)(1 +2 \alpha X + 3 X^2) $\\
\hline
\ttoo&2&$T_2$&$ (1- X)(1-2 X)(1 -4 \alpha X + 32 X^2)$\\
&&$T_3$&$ (1- X)(1-3 X)(1 +18 \alpha X + 243 X^2) $\\
\hline
\foot&2&$T_2$&$ (1- 2 X)(1-4 X)(1 -\beta X + 8 X^2)$\\
&&$T_3$&$ (1- 3 X)(1-9 X)(1 + (2\beta +6)X + 27 X^2) $\\
\hline
\multicolumn{4}{|l|}{$\dim S_2(31) = 2$, $ \dim S_4(31) = 2+5$}\\
\hline
\end{tabular}
\end{center}


{Here $\alpha$ satisfies $\alpha^2-\alpha-1=0$, and $\beta$ satisfies $\beta^2+5\beta+2=0$.}

\testentry
%
%
%
%
\begin{center}
\begin{tabular}{|p{40pt}|p{40pt}|p{40pt}|p{200pt}|}
\hline
\multicolumn{4}{|l|}{Level \textbf{33}. $\rrrank = 10$.}\\
\hline
\hline
\tott&1&$T_2$&$ (1-4 X)(1-8 X)(1 - X + 2 X^2)$\\
\hline
\ttoo&1&$T_2$&$ (1- X)(1-2 X)(1 -4 X + 32 X^2)$\\
\hline
\foot&1&$T_2$&$ (1- 2 X)(1-4 X)(1 +5 X + 8 X^2)$\\
\hline
\foot&1&$T_2$&$ (1- 2 X)(1-4 X)(1 + X + 8 X^2)$\\
\hline
\tott&3&$T_2$&$ (1-4 X)(1-8 X)(1 +2 X + 2 X^2)$\\
\hline
\ttoo&3&$T_2$&$ (1- X)(1-2 X)(1 +8 X + 32 X^2)$\\
\hline
\multicolumn{4}{|l|}{$\dim S_2(33) = 1$, $ \dim S_4(33) = 1+1+2+2$}\\
\hline
\end{tabular}
\end{center}


{The three dimensional eigenspaces are lifts of the weight two newform from level $11$.}

\testentry
%
%
%
%
\begin{center}
\begin{tabular}{|p{40pt}|p{40pt}|p{40pt}|p{200pt}|}
\hline
\multicolumn{4}{|l|}{Level \textbf{35}. $\rrrank = 7$.}\\
\hline
\hline
\tott&1&$T_2$&$ (1-4 X)(1-8 X)(1  + 2 X^2)$\\
\hline
\ttoo&1&$T_2$&$ (1- X)(1-2 X)(1  + 32 X^2)$\\
\hline
\tott&2&$T_2$&$ (1-4 X)(1-8 X)(1 -\alpha X + 2 X^2)$\\
\hline
\ttoo&2&$T_2$&$ (1- X)(1-2 X)(1 -4 \alpha X + 32 X^2)$\\
\hline
\foot&1&$T_2$&$ (1- 2 X)(1-4 X)(1 - X + 8 X^2)$\\
\hline
\multicolumn{4}{|l|}{$\dim S_2(35) = 1+2$, $ \dim S_4(35) = 1+2+3$}\\
\hline
\end{tabular}
\end{center}


{Here $\alpha$ satisfies $\alpha^2+\alpha-4=0$.}

\testentry
%
%
%
%
\begin{center}
\begin{tabular}{|p{40pt}|p{40pt}|p{40pt}|p{200pt}|}
\hline
\multicolumn{4}{|l|}{Level \textbf{37}. $\rrrank = 8$.}\\
\hline
\hline
\tott&1&$T_2$&$ (1-4 X)(1-8 X)(1 +2 X + 2 X^2)$\\
\hline
\ttoo&1&$T_2$&$ (1- X)(1-2 X)(1 +8 X + 32 X^2)$\\
\hline
\tott&1&$T_2$&$ (1-4 X)(1-8 X)(1  + 2 X^2)$\\
\hline
\ttoo&1&$T_2$&$ (1- X)(1-2 X)(1  + 32 X^2)$\\
\hline
\foot&4&$T_2$&$ (1- 2 X)(1-4 X)(1 -\alpha X + 8 X^2)$\\
\hline
\multicolumn{4}{|l|}{$\dim S_2(37) = 1+1$, $ \dim S_4(37) = 4+5$}\\
\hline
\end{tabular}
\end{center}


{Here $\alpha$ satisfies $\alpha^4+6\alpha^3-\alpha^2-16\alpha+6=0$.}

\testentry
%
%
%
%
\begin{center}
\begin{tabular}{|p{40pt}|p{40pt}|p{40pt}|p{200pt}|}
\hline
\multicolumn{4}{|l|}{Level \textbf{39}. $\rrrank = 10$.}\\
\hline
\hline
\tott&1&$T_2$&$ (1-4 X)(1-8 X)(1 - X + 2 X^2)$\\
\hline
\ttoo&1&$T_2$&$ (1- X)(1-2 X)(1 -4 X + 32 X^2)$\\
\hline
\tott&2&$T_2$&$ (1-4 X)(1-8 X)(1 -\alpha X + 2 X^2)$\\
\hline
\ttoo&2&$T_2$&$ (1- X)(1-2 X)(1 -4 \alpha X + 32 X^2)$\\
\hline
\foot&1&$T_2$&$ (1- 2 X)(1-4 X)(1  + 8 X^2)$\\
\hline
\foot&3&$T_2$&$ (1- 2 X)(1-4 X)(1 +5 X + 8 X^2)$\\
\hline
\multicolumn{4}{|l|}{$\dim S_2(39) = 1+2$, $ \dim S_4(39) = 1+2+2$}\\
\hline
\end{tabular}
\end{center}


{Here $\alpha$ satisfies $\alpha^2+2\alpha-6=0$.  The three
dimensional eigenspaces are lifts of the weight four newform from level $13$.}

\testentry
%
%
%
%
\begin{center}
\begin{tabular}{|p{40pt}|p{40pt}|p{40pt}|p{200pt}|}
\hline
\multicolumn{4}{|l|}{Level \textbf{41}. $\rrrank = 9$.}\\
\hline
\hline
\tott&3&$T_2$&$ (1-4 X)(1-8 X)(1 -\alpha X + 2 X^2)$\\
\hline
\ttoo&3&$T_2$&$ (1- X)(1-2 X)(1 -4 \alpha X + 32 X^2)$\\
\hline
\foot&3&$T_2$&$ (1- 2 X)(1-4 X)(1 -\beta X + 8 X^2)$\\
\hline
\multicolumn{4}{|l|}{$\dim S_2(41) = 3$, $ \dim S_4(41) = 3+7$}\\
\hline
\end{tabular}
\end{center}


{Here $\alpha$ satisfies $\alpha^3+\alpha^2-5\alpha-1=0$, and $\beta$ satisfies $\beta^3+3\beta^2-5\beta-3=0$.}

\testentry
%
%
%
%
\begin{center}
\begin{tabular}{|p{40pt}|p{40pt}|p{40pt}|p{200pt}|}
\hline
\multicolumn{4}{|l|}{Level \textbf{43}. $\rrrank = 10$.}\\
\hline
\hline
\tott&1&$T_2$&$ (1-4 X)(1-8 X)(1 +2 X + 2 X^2)$\\
\hline
\ttoo&1&$T_2$&$ (1- X)(1-2 X)(1 +8 X + 32 X^2)$\\
\hline
\tott&2&$T_2$&$ (1-4 X)(1-8 X)(1 -\alpha X + 2 X^2)$\\
\hline
\ttoo&2&$T_2$&$ (1- X)(1-2 X)(1 -4 \alpha X + 32 X^2)$\\
\hline
\foot&4&$T_2$&$ (1- 2 X)(1-4 X)(1 -\beta X + 8 X^2)$\\
\hline
\multicolumn{4}{|l|}{$\dim S_2(43) = 1+2$, $ \dim S_4(43) = 4+6$}\\
\hline
\end{tabular}
\end{center}


{Here $\alpha$ satisfies $\alpha^2-2=0$, and $\beta$ satisfies $\beta^4+4\beta^3-9\beta^2-14\beta+2=0$.}

\testentry
%
%
%
%
\begin{center}
\begin{tabular}{|p{40pt}|p{40pt}|p{40pt}|p{200pt}|}
\hline
\multicolumn{4}{|l|}{Level \textbf{47}. $\rrrank = 11$.}\\
\hline
\hline
\tott&4&$T_2$&$ (1-4 X)(1-8 X)(1 -\alpha X + 2 X^2)$\\
\hline
\ttoo&4&$T_2$&$ (1- X)(1-2 X)(1 -4 \alpha X + 32 X^2)$\\
\hline
\tott&3&$T_2$&$ (1-4 X)(1-8 X)(1 -\beta X + 2 X^2)$\\
\hline
\multicolumn{4}{|l|}{$\dim S_2(47) = 4$, $ \dim S_4(47) = 3+8$}\\
\hline
\end{tabular}
\end{center}


{Here $\alpha$ satisfies $\alpha^4-\alpha^3-5\alpha^2+5\alpha-1=0$, and $\beta$ satisfies $\beta^3+5\beta^2-2\beta-12=0$.}

\testentry 

\section{Interpretation of the numerical results}\label{interpretation}
\subsection{}
The first step, for each Hecke eigenvector $\beta \in H_5(\Gamma,\C)$
and for each prime $l$, is to write down and factor the Hecke
polynomial $P_l(X)$.

We then see that for the data computed so far, $\beta$ has one of the
following Galois representations attached.  As above, $p$ is a prime
not dividing $N$ or any of the $l$'s we are looking at.  

Let $k=2$ or $4$ and consider the continuous semisimple representation
$\sigma_k : G_{\Q}\rightarrow \GL_{2} (\Q_p)$ unramified outside $pN$
attached to a classical Hecke eigenform $f$ of weight $k$ and level $N'$
dividing $N$.  If $N$ is prime, we also assume $f$ has trivial
nebentypus.  Let $\rho = \epsilon^a \sigma_k \oplus \epsilon^b \oplus
\epsilon^c$ where $(k,a,b,c) = (2,0,2,3)$ or $(2,2,0,1)$, or
$(4,0,1,2)$.  If $f$ is an Eisenstein series, then $\rho = \chi_0\ \oplus
\chi_1 \epsilon \oplus \chi_2 \epsilon^2 \oplus \chi_3 \epsilon^3$, where the
$\chi_i$ are Dirichlet characters of conductor dividing $N$ with values in
a finite extension of $\F$, at least two of which are trivial.
  
Then for any $\beta$ there is some choice of such $\rho$ which is
apparently attached to $\beta$, in the sense that the Hecke polynomial
at $l$ equals the characteristic polynomial of $\Frob_l$ for all $l$
for which we computed the Hecke eigenvalues.

Thus it appears that none of our computed classes so far is cuspidal.
Therefore we should be able to related them to cohomology of the
boundary, either geometrically or in terms of Eisenstein series.  We
cannot do this thoroughly, because neither the cohomology of the
Borel-Serre boundary nor the theory of Eisenstein cohomology has been
sufficiently worked out for $\GL_{4}/\Q$.  This is not an easy task.
We can give the following indications.

From results of Moeglin-Waldspurger \cite{waldspurger}, we don't
expect any of our classes to be residues of Eisenstein series.  In the
framework of \cite{franke.schwermer} we can guess that our classes lie
either in the part of the cohomology indexed by the associate class of
parabolic subgroups of $\GL_{4}$ of type $(2,1,1)$ or in the part
indexed by the Borel subgroup.  Here the cuspidal data on the
$\GL_{2}$-factor of the Levi component of the first parabolic comes
from the appropriate classical cuspform of weight 2 or 4,  and we use
the appropriate power of the determinant on the $\GL_{1}$ factors.

\subsection{}
Geometrically, we make the following comments.  Let $M$ be the
quotient of the symmetric space for $\SL_{4} (\R )$ by $\Gamma$ and
let $\partial M$ be the boundary of its Borel-Serre compactification.
The covering of $\partial M$ by its faces gives a spectral sequence
for its cohomology.  The $E_2$ page has for its $(i,j)$-th term
$H^i(\text{Tits building}/\Gamma, H^j(\text{Fiber}))$.  An element of
that is an assignment: to every face $e'(P)$ of codimension $i$ we
assign an element of the cohomology in degree $j$ of $P \cap \Gamma$.
These assignments when restricted to a common face of codimension
$i+1$ must add up to 0.  Such an assignment gives a class in
$E_2^{i,j} (\partial M)$.  If it persists in the spectral sequence to
$E_\infty$, it will contribute to the cohomology $H^{i+j}(\partial
M,\C)$.  There remains the question as to whether this contribution is
the restriction of a class in $H^{i+j}(M,\C)$.

Note that if $P = LU$ is a Levi decomposition of $P$ then the spectral
sequence of the fibration for the cohomology of $P \cap \Gamma$
corresponding to this decomposition is known to degenerate at $E_2^{p,q} =
H^p(\Gamma_L, H^q(U \cap \Gamma)$ where
$\Gamma_L$ is the projection of $\Gamma$ to $L$.

The classes we have computed so far we expect to be coming in this way
from $\partial M$.  From the shape of the apparently associated Galois
representations, here is what we believe is their origin.  We only
sketch the constructions, since a detailed description would require a
thorough investigation of the cohomology of $\partial M$ for arbitrary
congruence subgroups of $\SL_{4}(\z)$.  First assume $f$ is a cuspform.

\subsection{}
The case where $\sigma$ has weight 2: By the Eichler-Shimura theorem,
the cuspform $f$ that has $\sigma$ attached shows up as a class
$\alpha$ in $H^1(\Delta,\C)$, where $\Delta$ is the classical
$\Gamma_0(N') \subset \SL_{2} (\Z )$.  First suppose $N'=N$.  Consider
the following element of $E_2^{0,5}(\partial M)$: On the standard
parabolic subgroup of type $(2,2)$ which is the stabilizer of the span
of $(e_1,e_2)$ in 4-space, we put the cohomology class $\alpha \times
1$ on the Levi component where we view the trivial coefficients of
$\alpha$ as the module $H^4(U \cap \Gamma)$.  One sees that there is a
unique class in the appropriate face corresponding to a $(3,1)$-type
parabolic subgroup that has the same restriction to the type $(2,1,1)$
face they have in common and restricts to 0 on the other faces.  Hence
these two glue up to give a class in the boundary.  The same
construction with the transposed parabolic subgroup also gives a
class, and these seem to account for all our $\beta$'s falling under
this case.  If $N' \ne N$, we choose a $(2,2)$- parabolic subgroup $P$
such that $\Gamma_L$ has level $N'$.  Then we imitate the construction
above.

\subsection{}
The case where $\sigma$ has weight 4: Here not every class we
construct on the boundary seems to lift to $H^5 (M)$, but only some of
them.  We don't know the reason for this.  The construction here
creates a class in $E_2^{1,4}(\partial M)$, and if this correctly
describes what we have computed, our computed classes of this type
would be ghost classes.  That is, the corresponding cohomology class
in $H^5(M)$ restricts nontrivially to the boundary of $M$, but it
restricts to 0 on each face of the boundary, since it is coming from a
class in $E_2^{i,j}(\partial M)$ with $i > 0 $.

For this construction, one first chooses a parabolic subgroup $P$ of
type $(2,1,1)$.  Note that $H^3((U \cap \Gamma),\C)$ contains a
$\Gamma_L$-submodule isomorphic to $\mathcal V_2$, the homogeneous complex
polynomials of degree 2 on 2-space, after we identify $\Gamma_L$ with
a subgroup of $\GL_{2} (\Z ) \times \GL_{1} (\Z )\times \GL_{1} (\Z
)$.  By the Eichler-Shimura theorem, the weight $4$ cuspform $f$ that
has $\sigma$ attached shows up in $H^1(\Delta,\mathcal V_2)$.  Thus we can view
it as in $H^1(\Gamma_L, H^3(U \cap \Gamma))$.

We can do this on three $P$'s which are not conjugate to each other in
such a way that they, together with certain classes on $(3,1)$-type
parabolic subgroups, all glue together to give a class in
$E_2^{1,4}(\partial M)$.  The details are left for the reader.

\subsection{}
Finally, the case where $\sigma$ is an Eisenstein series is
harder to understand.  We haven't worked out exactly how the gluing process
goes in this case, so we're not sure what stratum of the spectral sequence
is occupied by the corresponding boundary classes.  

We note that in every case, $\rho$ restricted to an inertia
subgroup of $G_\Q$ at $p$ has the form $1 \oplus \epsilon \oplus \epsilon^2
\oplus \epsilon^3$, which is consistent with the conjecture of Ash-Sinnott
\cite{ash.sinnott}, since the coefficient module of the cohomology
classes we consider is the trivial module. 


\providecommand{\bysame}{\leavevmode\hbox to3em{\hrulefill}\thinspace}

\end{document}